\documentclass[bj,authoryear]{imsart}

\RequirePackage{amsthm,amsmath,amsfonts,amssymb}
\allowdisplaybreaks
\RequirePackage{graphicx}%% uncomment this for including figures

\usepackage{subcaption}
\usepackage{float} % [htbp] option for figure

\usepackage{caption}
\usepackage{booktabs}
\usepackage{float}
\usepackage{multirow}
\usepackage{diagbox}
\usepackage{bbm}

\startlocaldefs
\theoremstyle{plain}

\newtheorem{theorem}{Theorem}[section]
\newtheorem{lemma}[theorem]{Lemma}
\newtheorem{Corollary}[theorem]{Corollary}
\newtheorem{Proposition}[theorem]{Proposition}

\newtheorem{Assumption}[theorem]{Assumption}
\theoremstyle{remark}

\newtheorem{remark}{Remark}
\newcommand{\bmu}{\boldsymbol{\mu}}
\newcommand{\bSigma}{\boldsymbol{\Sigma}}
\newcommand{\bPhi}{\boldsymbol{\Phi}}
\newcommand{\E}{\mathbb{E}}
\newcommand{\bA}{\mathbf{A}}
\newcommand{\bB}{\mathbf{B}}
\newcommand{\bC}{\mathbf{C}}
\newcommand{\bD}{\mathbf{D}}

\newcommand{\bF}{\mathbf{F}}

\newcommand{\bH}{\mathbf{H}}
\newcommand{\bI}{\mathbf{I}}

\newcommand{\bM}{\mathbf{M}}

\newcommand{\bQ}{\mathbf{Q}}
\newcommand{\bR}{\mathbf{R}}
\newcommand{\bS}{\mathbf{S}}

\newcommand{\bX}{\mathbf{X}}
\newcommand{\bY}{\mathbf{Y}}

\newcommand{\be}{\mathbf{e}}

\newcommand{\br}{\mathbf{r}}

\newcommand{\bx}{\mathbf{x}}
\newcommand{\by}{\mathbf{y}}

\newcommand{\calN}{\mathcal{N}}

\newcommand{\Prob}{\mathbb{P}}
\newcommand{\Expe}{\mathbb{E}}
\newcommand{\Var}{\mathrm{Var}}

\newcommand{\tr}{\mathrm{tr}}

\newcommand*{\dif}{\mathop{}\!\mathrm{d}}
\renewcommand{\bar}{\overline}
\usepackage{dsfont}

\endlocaldefs

\begin{document}

\begin{frontmatter}
\title{On eigenvalues of a  renormalized sample correlation matrix}
%\title{A sample article title with some additional note\thanksref{T1}}
\runtitle{On eigenvalues of a  renormalized sample correlation matrix}
%\thankstext{T1}{A sample of additional note to the title.}

\begin{aug}
\author[A]{\fnms{Qianqian}~\snm{Jiang}* \ead[label=e1]{jqq172515@gmail.com}}
\author[A]{\fnms{Junpeng}~\snm{Zhu}* \ead[label=e2]{zhujp@sustech.edu.cn}}
\author[A]{\fnms{Zeng}~\snm{Li}\ead[label=e3]{liz9@sustech.edu.cn}}

%\footnotetext{* These authors contributed equally.}

%%%%%%%%%%%%%%%%%%%%%%%%%%%%%%%%%%%%%%%%%%%%%%
%% Addresses                                %%
%%%%%%%%%%%%%%%%%%%%%%%%%%%%%%%%%%%%%%%%%%%%%%

\address[A]{Department of Statistics and Data Science, Southern University of Science and Technology \footnote{These authors contributed equally}
\printead[presep={,\ }]{e1,e2,e3}}

\end{aug}

\begin{abstract}
This paper studies the asymptotic spectral properties of a renormalized sample correlation matrix, including the limiting spectral distribution, the properties of largest eigenvalues, and the central limit theorem for linear spectral statistics. All asymptotic results are derived under a unified framework
where the dimension-to-sample size ratio  $p / n \rightarrow c \in(0, \infty]$.
Based on our CLT result, we propose an independence test statistic capable of operating effectively in both high and ultrahigh dimensional scenarios.
Simulation experiments demonstrate the accuracy of theoretical results.
\end{abstract}



\end{frontmatter}

%%%%%%%%%%%%%%%%%%%%%%%%%%%%%%%%%%%%%%%%%%%%%%
%%%% Main text entry area:

%\tableofcontents

\section{Introduction}
 Let us consider the widely used independent components (IC) model for the population $\mathbf{x}$, admitting the following stochastic representation
%$$
%\mathbf{x}=\bmu+{\bSigma}^{\frac{1}{2}} \mathbf{z},
%$$
$$
\mathbf{y}=\bmu+{\bSigma}^{\frac{1}{2}} \mathbf{x},
$$
where $\bmu \in \mathbb{R}^p$ denotes the population mean and $\mathbf{x} \in \mathbb{R}^p$ is a  random vector with independent and identically distributed (i.i.d.) components with zero mean and unit variance. Let $\mathbf{y}_1, \ldots, \mathbf{y}_n$ be $n$ i.i.d. observations from this population and $\bY=(\mathbf{y}_1, \ldots, \mathbf{y}_n)$ be the $p\times n$ data matrix. The sample correlation matrix $\bR_n$ can be written as 
%$$
%\bR_n=\bY_n^{*}\bY_n, 
%$$
%where 
%$$
%\mathbf{Y}_n=\left(\frac{\mathbf{y}_1-\overline{\mathbf{y}}_1}{\left\|\mathbf{y}_1-\overline{\mathbf{y}}_1\right\|}, \frac{\mathbf{y}_2-\overline{\mathbf{y}}_2}{\left\|\mathbf{y}_2-\overline{\mathbf{y}}_2\right\|}, \ldots, \frac{\mathbf{y}_p-\overline{\mathbf{y}}_p}{\left\|\mathbf{y}_p-\overline{\mathbf{y}}_p\right\|}\right),~\mathbf{y}_i=\be_i^{\top}\bX,~ \bar\by_i=\frac{1}{n}\mathbf{1}_n\mathbf{1}_n^{\top}\by_i,
%$$
%$$
%\mathbf{Y}_n=\left(\frac{\bPhi\mathbf{y}_1}{\left\|\bPhi\mathbf{y}_1\right\|}, \frac{\bPhi\mathbf{y}_2}{\left\|\bPhi\mathbf{y}_2\right\|}, \ldots, \frac{\bPhi\mathbf{y}_p}{\left\|\bPhi\mathbf{y}_p\right\|}\right),~\bPhi=\mathbf{I}_n-\frac{1}{n}\mathbf{1}_n \mathbf{1}_n^{\top},
%~\mathbf{y}_i^*=\be_i^{\top}\bX,~ i\in\{1,\ldots,p\}.
%$$

$$
\bR_n=\bD_n^{\frac{1}{2}}\bS_n\bD_n^{\frac{1}{2}},
$$
where
$$
\bD_n^{\frac{1}{2}}=\operatorname{Diag}\left(\frac{1}{\sqrt{s_{11}}}, \frac{1}{\sqrt{s_{22}}}, \ldots, \frac{1}{\sqrt{s_{p p}}}\right),~\bS_n=\frac{1}{N}\bY\bPhi\bY^{\top},~
\boldsymbol{\Phi}=\mathbf{I}_n-\frac{1}{n} \mathbf{1}_n \mathbf{1}_n^{\top},~N=n-1.
$$
Here 
$s_{kk}=\be_k^{\top}\bS_n\be_k,~k=1, \ldots, p$, 
$\mathbf{e}_i \in \mathbb{R}^p$ denotes the vector with the $i$ the element being 1 and all others being 0, and $\mathbf{1}_n=(1, \ldots, 1)^{\top}$ in $\mathbb{R}^n$.

The eigenvalues of $\bR_n$,
$\lambda^{\bR_n}_1\geq\dots\geq\lambda^{\bR_n}_p$,
serve as important statistics and often play crucial roles in the inference on population correlation matrix $\mathcal{R}$, see \cite{Anderson}.  Consider the following regime,
\begin{align}\label{mp}
 n\to \infty, p=p_n\to \infty, p/n\to c\in(0,\infty),
\end{align}
referred to as the Marčenko-Pastur (MP) regime. 
For $\mathcal{R}=\bI_p$,
\cite{jiang2004limiting}
demonstrated that
the empirical spectral distribution (ESD) of $\bR_n$, $F^{\mathbf{R}_n}(x)=\frac{1}{p} \sum_{i=1}^p \mathbbm{1}_{\left\{\lambda_i\left(\mathbf{R}_n\right) \leq x\right\}}$, converges weakly to the Marchenko-Pastur (MP) law with probability one. 
The extreme eigenvalues of $\bR_n$ were studied in  \cite{MR2591901} and 
\cite{MR2988403}. 
Additionally, \cite{gao2017high} established the central limit theorem (CLT) for the linear spectral statistics (LSS) of $\bR_n$, i.e., $\int f(x) \mathrm{d} F^{\mathbf{R}_n}(x)= \sum_{i=1}^p f(\lambda_i^{\bR_n})/p$ where $f(\cdot)$ is a smooth function. For a general $\mathcal{R}$, the limiting spectral distribution (LSD) of $\bR_n$, the limit of ESD, can be found in \cite{Karoui2009} and 
the CLT for LSS was studied in \cite{mestre2017correlation,jiang2019determinant,yin2023central,MR4478187}.
All these studies are conducted under the MP regime \eqref{mp}, i.e., $p / n \rightarrow c \in(0, \infty)$.

However, in the ultrahigh dimensional case where $p \gg n$, the eigenvalues of $\mathbf{R}_n$ exhibit behaviors markedly different from those in the MP regime. Properties of eigenvalues of sample correlation matrix when $p \gg n$ remain largely unknown in current literature. 
Existing studies on eigenvalue behavior of ultrahigh dimensional matrix focus on 
sample covariance matrix, see  \cite{bai1988convergence,WANG201425,chen2015clt,bao2015,Qiuclt}. These works heavily rely on the linear independent component structure and zero mean assumption $\bmu=0$ which suggest that the renormalized sample covariance matrix 
$\tilde{\bS}_n=\sqrt{\frac{p}{n}}\left(\frac{1}{p}\mathbf{Y}^{\top}_0 \mathbf{Y}_0 -\bI_n\right)$, $\mathbf{Y}_0=\bY-\bmu\mathbf{1}_n^{\top}$
shares many spectral properties with Wigner matrix. In contrast, due to the nonlinear dependence introduced by the normalization inherent in the sample correlation matrix and the presence of a nonzero population mean, the techniques and results developed for ultrahigh dimensional covariance matrices cannot be directly extended to the correlation matrix. To fill this gap, we consider the sample correlation matrix under a new regime where $p / n \rightarrow \infty$ as $n \rightarrow \infty$. In this scenario, unlike the MP regime, most eigenvalues of the matrix $\bR_n$ are zero, and all non-zero eigenvalues diverge to infinity. To address this, we renormalize the sample correlation matrix as follows:
 $$
\bB_n=\sqrt{\frac{p}{N}}\left(\frac{1}{p} \bPhi\bY^{\top}\bD_n\bY\bPhi-\boldsymbol{\Phi}\right).
$$
$\mathbf{B}_n$ is $n \times n$ and has $n-1$ non-zero eigenvalues, which connect to the non-zero eigenvalues of $\mathbf{R}_n$ through the following identity:

$$
\lambda^{\mathbf{B}_n}=\sqrt{\frac{N}{p}} \lambda^{\mathbf{R}_n}-\sqrt{\frac{p}{N}}.
$$

This paper investigates the eigenvalues of the renormalized random matrix $\mathbf{B}_n$ when $\mathcal{R}=\bI_p$, allowing for the dimension $p$ to be comparable to or much larger than the sample size $n$ , such that
$$
n \rightarrow \infty, p=p_n \rightarrow \infty, p / n \rightarrow c \in(0, \infty].
$$
Firstly, we propose a unified LSD of $\bB_n$ in both $p/n\to c\in(0,\infty)$ and $p/n\to \infty$. 
Secondly, we studied the properties of 
$\lambda_1^{\bB_n}$,
the largest eigenvalue of $\bB_n$.
Thirdly, we establish CLT for LSS of $\bB_n$ under the unified framework, which covers the results in \cite{gao2017high} as a special case. Last but not least, our theoretical findings are further applied to the independence test for both high and ultrahigh dimensional random vectors. Specifically, we propose a test statistic that remains effective when  $p / n \rightarrow c \in(0, \infty]$. 

  % The theoretical analysis of the renormalized sample correlation matrix $\bB_n$ in ultrahigh dimensional settings presents significant challenges. This complexity stems from the nonlinear dependence structure introduced by the normalization process, making the study of this random matrix more intricate, even when $\mathcal{R} = \bI_p$. 

In this paper, our primary contribution 
 is to establish the asymptotic theory for eigenvalues of the renormalized sample correlation matrix $\bB_n$ when $p/n\to \infty$. In addition, we provide a unified representation of the limiting results that hold for both $p/n\to \infty$ and $p/n\to c\in(0,\infty)$. Theoretical analysis of  $\bB_n$ in ultrahigh dimensional settings presents significant challenges due to the nonlinear dependence structure introduced by the normalization process, which makes the study of this random matrix more intricate, even when \(\mathcal{R} = \bI_p\). 
Under the MP regime \eqref{mp}, \cite{jiang2004limiting,Karoui2009,heiny2022large,gao2017high} showed that the correlation matrix $\bR_n$ share the same LSD and properties of the largest eigenvalue as the sample covariance $\bS_n$, by using $\|\bD_n-\bI_p\|\|\bS_n\|$ to control the difference between the sample correlation matrix $\bR_n$ and the sample covariance matrix $\bS_n$. However, in the ultrahigh dimensional setting, since $\|\bS_n\|$ tends to infinity, this approach becomes ineffective. Instead, we investigate the convergence of Stieltjes transform of ESD of $\bB_n$ to obtain the LSD.  In addition, we require a unified moment assumption to control the probability that the largest eigenvalue $\lambda_1^{\bB_n} $ lies outside the support of LSD. Moreover, when $p/n\to c\in(0,\infty)$, \cite{Pan2008central} used $c_n=p/n$ to characterize the CLT for LSS while \cite{yin2023central,MR4478187}
used $c_N=p/N$. In fact, they are equivalent because, in the high-dimensional setting \eqref{mp}, $c_n-c_N=O(1 / n)$. However, when $p/n\to \infty$, $c_n-c_N=p /(n N)$ may diverge to infinity. Therefore, we must handle $c_n$ and $c_N$ with extra caution 
and we derive a novel determinant equivalent form
for the resolvent of the renormalized correlation matrix 
when $p/n\to \infty$.

The rest of the paper is organized as follows. Section 2 details our main results, including unified 
LSD, properties of the largest eigenvalue and 
CLT for LSS. Section 3 discusses the application of our CLT to independence test. Section 4 presents simulations. Technical proofs are detailed in Section 5 and 
the Supplementary Material.

\section{Main Results}

\subsection{Preliminaries}
%For any $n \times n$ Hermitian matrix $\bB_n$ with eigenvalues $\lambda_1,\ldots,\lambda_n$, its \emph{empirical spectral distribution} (ESD) is defined by
%\begin{equation*}
	%F^{\bB_n}(x)=\frac{1}{n}\sum^{n}_{i=1}I_{\{\lambda_i (\bB_n)\leq   x\}},
%\end{equation*}
%where $I_{\{\cdot\}}$ denotes the indicator function. If $F^{\bB_n}(x)$ converges to a non-random limit $F(x)$ as $n\rightarrow\infty$, we call $F (x)$ the \emph{limiting spectral distribution} of $\bB_n$. The LSD of $\bB_n$ is described in terms of its Stieltjes transform. 

For any measure $G$ supported on the real line, the Stieltjes transform of $G$ is defined as
\[
s_G(z)=\int\frac{1}{x-z}\dif G(x),\quad z \in \mathbb{C}^+,
\]
where $\mathbb{C}^{+}=\{z \in \mathbb{C}: \Im(z)>0\}$ denotes the upper complex plane.

As for the LSD of $\bR_n$ with $\mathcal{R}=\bI_p$ when $p/n\to c\in(0,\infty)$, 
\cite{jiang2004limiting} showed the ESD of $\bR_n$ converges with probability $1$ to the Marčenko-Pastur law $F_{MP}(x)$, whose density function has an explicit expression
$$
f_{MP}(x)= \begin{cases}\frac{1}{2 \pi x c} \sqrt{(b-x)(x-a)} & a \leqslant x \leqslant b; \\ 0 & \text { otherwise, }\end{cases}
$$
and a point mass $1-1 / c$ at the origin if $c>1$, where $a=(1-\sqrt{c})^2$ and $b=(1+\sqrt{c})^2$.
And the Stieltjes transform of $F_{MP}(x)$ is 
\begin{align}\label{mp-lsd}
s_{MP}(z)=\frac{c-1-z+\sqrt{(z-c-1)^2-4 c}}{2cz}+\frac{1-c}{cz},\quad 
z \in \mathbb{C}^{+}.
\end{align}

%Many classes of statistics related to the eigenvalues of $\bB_{n}$ are important for multivariate inference, particularly functionals of the ESD. To explore this, for any function $f$ defined on $[0,\infty)$, we consider the \emph{linear spectral statistics} of $\bB_{n}$ given by
%\[
%\int f(x)\dif F^{\bB_{n}}(x)=\frac{1}{n}\sum_{i=1}^n f\bigl(\lambda_i\bigr),
%\]
%where $\lambda_i$, $i=1,\ldots,n$, are eigenvalues of $\bB_{n}$.

%In this paper, we study the asymptotic spectral properties of $\bB_{n}$, including the LSD (see, Theorem \ref{thm:LSD}), the behavior of extreme eigenvalues (see, Proposition \ref{cor-extreigen}), and the CLT for LSS (see, Theorem \ref{cltmain}). 

\subsection{LSD of $\bB_n$}
In this section, we provide a unified LSD of the renormalized 
sample correlation matrix
$\bB_n$ when $p/n\to c\in(0,\infty]$.

%Before diving into linear functionals of eigenvalues of $\bB_{n}$, we first explore its LSD and extreme eigenvalues. Specifically, suppose the following assumptions hold.
\begin{Assumption}\label{asp1}
  Let $\mathbf{X}=\left(\mathbf{x}_{1}, \ldots, \mathbf{x}_{n}\right)_{p\times n}=\left(x_{i j}\right)$, which consists of $p \times n$ i.i.d. variables satisfying

$$
\E\left(x_{i j }\right)=0, \quad \E\left|x_{i j }\right|^2=1, \quad \E\left|x_{i j }\right|^4=\kappa<\infty.
$$

\begin{Assumption}\label{asp2}
   The population covariance matrix $\bSigma$ is diagonal.
\end{Assumption}

\end{Assumption}
\begin{Assumption}\label{asp3}
The dimension $p$ is function of sample size $n$ and both tend to infinity such that
   $$p/n\to c\in (0,\infty],~
    p \asymp n^t,~t\geq 1.$$
\end{Assumption}

\begin{theorem}\label{thm:LSD}
    Under Assumptions \ref{asp1} - \ref{asp3}, with probability one, the ESD of $\bB_{n}$ converges weakly to a (non-random) probability measure $F^c(x)$, which has a density function
    \begin{equation*}
     f^{c}(x) = \begin{cases}
		\frac{\sqrt{4 -x^2-c^{-1}+2 xc^{-1/2}}}{2 \pi(1+ xc^{-1/2})},&\text{if }x\in\left[\frac{1}{\sqrt{c}}-2,\frac{1}{\sqrt{c}}+2\right],\\
		0,& \text{otherwise},
	\end{cases}
\end{equation*}
and 
has a point mass $1-c$ at the point $-\sqrt{c}$ if $0<c\leq 1$.
The Stieltjes transform of $F^c(x)$ is
\begin{align}\label{s_c(z)}
    s_{c}(z) = \frac{-(z+c^{-1/2})+\sqrt{(z+2-c^{-1/2})(z-2-c^{-1/2})}}{2(1+c^{-1/2}z)},~ z \in \mathbb{C}^{+}.
\end{align}
Moreover, the expression of the moments are
$$\int_{-\infty}^{+\infty}x^k f^{c}(x)\dif x=\sum_{s=0}^k(-1)^s\binom{k}{s}c^{-k/2+s+1}\beta_{k-s}+(1-c)(-\sqrt{c})^k,\quad k\geq1,$$
where $\beta_0=1$  and $\beta_j=\sum_{r=0}^{j-1} \frac{1}{r+1}\binom{j}{r}\binom{j-1}{r} c^r$ for $j\geq1$.
% are moments computed by continous density function in standard MP law.
 %, and $f^{c_N}(x)$ substitutes $c_N$ for $c$ in $f^{c}(x)$.
\end{theorem}
\begin{remark}
    Theorem \ref{thm:LSD} provides a unified LSD of $\bB_n$ when $p/n\to c\in (0,\infty]$. This result is consistent with the  MP law of $\bR_n$ when $p/n\to c\in (0,\infty)$ in  \eqref{mp-lsd}. %Since $f^c(x)$ converges to the semicicle law uniformly on a bounded interval as $c\to\infty$, the moments of $f^c(x)$ converge to moments of semicicle law.
\end{remark}
The following theorem shows the  result when $p/n\to \infty$, which, to the best of our knowledge, is presented here for the first time.
\begin{theorem}\label{ult-thm:LSD}
Under Assumptions \ref{asp1} - \ref{asp3} and $p \asymp n^t, t>1$, with probability one, the ESD of $\bB_{n}$ converges weakly to the semicircular law $F(x)$ with density function
\begin{equation}\label{eq:sc_density}
	f(x) = \begin{cases}
		\frac{1}{2\pi} \sqrt{4-x^2},&\text{if }x\in[-2,2],\\
		0,& \text{otherwise,}
	\end{cases}
	\end{equation}
 and Stieltjes transform $s(z)=\frac{-z+\sqrt{z^2-4}}{2},~z \in \mathbb{C}^{+}$. Moreover, the expression of the moments are $$\int_{-\infty}^\infty x^k \cdot \frac{1}{2 \pi} \sqrt{4-x^2} d x= \begin{cases}\frac{1}{k/2+1}\binom{k}{k/2}, & \text { if } k \text { is even, } \\ 0, & \text { if } k \text { is odd. }\end{cases}$$
    \end{theorem}

\iffalse
\begin{Proposition}\label{cor-extreigen}
    Under Assumptions \ref{asp1} and \ref{asp3}, 
 %    we have as $n\to\infty$
 %    \begin{equation}\label{eq:extreme-eig-CoDa}
	% 	\lambda_{\max}(\bB_{n}) \convas 2 \qquad \text{and}\qquad \lambda_{\min} (\bB_{n}) \convas -2,
	% \end{equation}
 % where $\lambda_{\max}(\bB_{n})$ is the largest eigenvalue of $\bB_{n}$, and $\lambda_{\min}(\bB_{n})$ is the smallest non-zero eigenvalue of $\bB_{n}$. 
for any $\ell>0$, $\eta_1>2$ and $\eta_2<-2$, we have
	\[
		\Prob\bigl(\lambda_{\max}(\bB_{n})\geq  \eta_1\bigr) = o(n^{-\ell})
	\qquad\text{and}\qquad
		\Prob\bigl(\lambda_{\min}(\bB_{n})\leq   \eta_2\bigr) = o(n^{-\ell}).
	\]
\end{Proposition}
\fi

\subsection{The largest eigenvalue of $\bB_n$}
In this section, we study the properties of $\lambda_1^{\bB_n}$, the largest eigenvalue of $\bB_n$, when $p \asymp n^t, t \geq 1$.

\renewcommand{\thetheorem}{2.1*} % Change the numbering format

\begin{Assumption}\label{asp1*}
  Let $\mathbf{X}=\left(\mathbf{x}_{1}, \ldots, \mathbf{x}_{n}\right)_{p\times n}=\left(x_{i j}\right)$, which consists of $p \times n$ i.i.d. variables satisfying

$$
\E\left(x_{i j }\right)=0, \quad \E\left|x_{i j }\right|^2=1, \quad \E\left|x_{i j }\right|^4=\kappa,\quad
\E\left|x_{i j }\right|^{2(t+1)}<\infty.
$$
\end{Assumption}
\setcounter{theorem}{5}
\renewcommand{\thetheorem}{\thesection.\arabic{theorem}}

\begin{remark}
 Compared with the assumptions in the literature \citep{gao2017high,yin2023central,MR4478187,jiang2004limiting,MR2988403} where $p / n \rightarrow c \in(0, \infty)$, Assumption \ref{asp1*} is not stronger. In fact, when $t=1$, the moment condition $\mathbb{E}\left|x_{i j}\right|^{2(t+1)}<\infty$ reduces to a finite fourth moment, which coincides with the standard assumption in random matrix theory .
    %Assumption \ref{asp1*} plays a crucial role in analyzing the properties of $\lambda_1^{\mathbf{B}_n}$. When $p/n \to c \in (0, \infty)$ with $t=1$, Assumption \ref{asp1*} reduces to Assumption \ref{asp1}, which is a standard assumption in random matrix theory \citep{BSbook}.
\end{remark}

\begin{theorem}\label{cor-extreigen}
    Under Assumptions \ref{asp1*},  \ref{asp2} and \ref{asp3}, we have
    \begin{itemize}
      \item [(i)] $\lambda_{1}\left(\mathbf{B}_n\right) \rightarrow 2+\frac{1}{\sqrt{c}} \quad$ a.s.;
      \item[(ii)]
      for any $\epsilon>0, \ell>0$,
if 
%\begin{align}*\label{con-max}
$\left|x_{i j}\right| \leq \delta_n (n p)^{1/(2t+2)},$ where $\delta_n \rightarrow 0, ~\delta_n (n p)^{1/(2t+2)}\rightarrow \infty$, as $ n \rightarrow \infty,$
%\end{align}
then 
$$\Prob\left(\lambda_{1}\left(\mathbf{B}_n\right) \geq 2+\frac{1}{\sqrt{c}}+\epsilon\right)=\mathrm{o}\left(n^{-\ell}\right).
$$
    \end{itemize}
\end{theorem}

\begin{remark}
Theorem \ref{cor-extreigen}  is consistence with the results of $\lambda_1^{\bR_n}$ when $p / n \rightarrow c \in(0, \infty)$ in Theorem 1.1 of \cite{jiang2004limiting} and Lemma 7 of 
\cite{gao2017high}. 
%The following Corollary shows the
%simplified result when $p / n \rightarrow \infty$.
\end{remark}

\iffalse
\begin{Corollary}
    Under Assumptions \ref{asp1*},  \ref{asp2}, \ref{asp3} and
    $p \asymp n^t, t>1$, we have
    \begin{itemize}
      \item [(i)] $\lambda_{1}\left(\mathbf{B}_n\right) \rightarrow 2\quad$ a.s.;
      \item[(ii)]
      for any $\epsilon>0, \ell>0$,
$
p\left(\lambda_{1}\left(\mathbf{B}_n\right) \geq 2+\epsilon\right)=\mathrm{o}\left(n^{-\ell}\right),
$
when \eqref{con-max} holds.
    \end{itemize}
\end{Corollary}
\fi

\subsection{CLT for LSS of $\bB_n$}
In this section, we focus on linear spectral statistic of $\bB_{n}$, i.e. $\frac{1}{n}\sum_{i=1}^nf(\lambda_i)$,
where $f$ is an analytic function on $[0, \infty)$. Since $F^{\bB_n}$ converges to $F^c$ almost surely, we have
$$\frac{1}{n}\sum_{i=1}^nf(\lambda_i)\to\int f(x)\dif F^c(x).$$ 
We explore second order fluctuation of $\frac{1}{n}\sum_{i=1}^nf(\lambda_i)$ describing how such LSS converges to its first order limit. Consider a renormalized functional:
$$G_{n}(f)=n\int_{-\infty}^{+\infty} f(x)\dif \left\{F^{\bB_{n}}(x)- F^{c_N}(x)\right\}+\frac{1}{2 \pi i}\oint_{\mathcal{C}} f(z)\Theta_n (s_{c_N}(z))\dif z,$$
where $F^{c_N}(x)$ and $s_{c_N}(z)$ serve as finite-sample proxies for $F^c(x)$ and $s_c(z)$ in \eqref{s_c(z)}, by substituting $c$ with $c_N=p/N$, \begin{align}
    \Theta_n (s_{c_N}(z)) &= 
    2 c_N^{-\frac{1}{2}} g_{c_N}^2(z) h_{c_N}(z) s^2_{c_N}(z)d_{c_N}(z)- c_N^{-\frac{1}{2}} g_{c_N}^2(z) h_{c_N}(z)  s_{c_N}^{\prime}(z)d_{c_N}(z)\nonumber\\
    &+\frac{1}{\sqrt{c_N}+z}g_{c_N}(z)h_{c_N}(z)d_{c_N}(z)- \frac{\sqrt{c_N}}{z(\sqrt{c_N}+z)}\nonumber\\
    &+\frac{n}{N}g_{c_N}(z)h_{c_N}(z) s_{c_N}(z)d_{c_N}(z)
    - \frac{c_N^{\frac32}}{-c_N^{-\frac{1}{2}}s_{c_N}(z)l^{-1}_{c_N}+\left(c_N+\sqrt{c_N} z\right)},
\end{align}
\iffalse
\begin{align}\label{real-theta}
          \Theta_n (s_{c_N}(z)) &= \frac{s^3_{c_N}(z)}{-1+s^2_{c_N}(z)} - \frac{\sqrt{c_N}}{z(\sqrt{c_N}+z)}-\frac{1}{\sqrt{c_N}}g^2_{c_N}(z)h_{c_N}(z)s^{\prime}_{c_N}(z)F_{c_N}(z)\nonumber\\
          &+\frac{n}{N}\frac{1}{\sqrt{c_N}+z}g_{c_N}(z)h_{c_N}(z)F_{c_N}(z), 
          % \ \ \text{under the real case},\nonumber
          % \Theta_n (s_{c_N}(z)) &= |\Psi|^2\frac{s^3_{c_N}(z)}{-1+s^2_{c_N}(z)} - \frac{\sqrt{c_N}}{z(\sqrt{c_N}+z)}\label{complex-theta}\\
          % &-\sqrt{c_N}|\Psi|^2\frac{s_{c_N}(z)h_{c_N}(z)}{h^{-2}_{c_N}(z)-c_N\Psi^2g^2_{c_N}(z)(c_N+\sqrt{c_N}z)^{-2}} F_{c_N}(z)\nonumber\\
          % &+\frac{n}{N}\frac{1}{\sqrt{c_N}+z}g_{c_N}(z)h_{c_N}(z)F_{c_N}(z),  \ \ \text{under the complex case},\nonumber
     \end{align}
     \fi
and 
\begin{align}
    h_{c_N}(z) &= \frac{1}{1+\frac{1}{\sqrt{c_N}}s_{c_N}(z)+\frac{1-c_n}{c_N+\sqrt{c_N}z}},\nonumber\\
     g_{c_N}(z) &= -\frac{c_N+\sqrt{c_N}z}{c_n}\left\{\frac{s_{c_N}(z)}{\sqrt{c_N}}+\frac{1-c_n}{c_N+\sqrt{c_N}z}\right\},\nonumber\\
     %F_{c_N}(z) &= \frac{\{1+\sqrt{c_N}U_{c_N}(z)s^{-1}_{c_N}(z)\}k^{-1}_{c_N}(z)l^{-1}_{c_N}(z)}{1+k^{-1}_{c_N}(z)l^{-1}_{c_N}(z)(c_N+\sqrt{c_N}z)}\nonumber\\
     d_{c_N}(z) &=- \frac{c_n h^{-1}_{c_N}(z)s_{c_N}(z)}{\sqrt{c_N}-l^{-1}_{c_N}(z)s_{c_N}(z)(c_N+\sqrt{c_N}z)},\nonumber\\
      l_{c_N}(z) &= \frac{h_{c_N}(z)}{c_n}\left[1+\frac{\sqrt{c_N}}{s_{c_N}(z)}\left\{c_n + \frac{c_n(1-c_n)}{c_N+\sqrt{c_N}z} + \left(c_n-1\right)\frac{s_{c_N}(z)}{\sqrt{c_N}}\right\}\right],\nonumber
     % U_{c_N}(z) &= c_n + \frac{c_n(1-c_n)}{c_N+\sqrt{c_N}z} + (c_n-1)\frac{s_{c_N}(z)}{\sqrt{c_N}}.\nonumber
\end{align}
% \begin{align}
%     h_{c_N}(z) &= \frac{1}{1+\frac{1}{\sqrt{c_N}}s_{c_N}(z)+\frac{1-c_n}{c_N+\sqrt{c_N}z}},\nonumber\\
%      g_{c_N}(z) &= -\frac{c_N+\sqrt{c_N}z}{c_n}\{\frac{s_{c_N}(z)}{\sqrt{c_N}}+\frac{1-c_n}{c_N+\sqrt{c_N}z}\},\nonumber\\
%      %F_{c_N}(z) &= \frac{\{1+\sqrt{c_N}U_{c_N}(z)s^{-1}_{c_N}(z)\}k^{-1}_{c_N}(z)l^{-1}_{c_N}(z)}{1+k^{-1}_{c_N}(z)l^{-1}_{c_N}(z)(c_N+\sqrt{c_N}z)}\nonumber\\
%      F_{c_N}(z) &= \frac{c_n h^{-1}_{c_N}(z)k^{-1}_{c_N}(z)}{1+k^{-1}_{c_N}(z)l^{-1}_{c_N}(z)(c_N+\sqrt{c_N}z)},\ k_{c_N}(z) = -\frac{\sqrt{c_N}}{s_{c_N}(z)},\nonumber\\
%       l_{c_N}(z) &= \frac{h_{c_N}(z)}{c_n}\{1+\frac{\sqrt{c_N}}{s_{c_N}(z)}U_{c_N}(z)\},\nonumber\\
%      U_{c_N}(z) &= c_n + \frac{c_n(1-c_n)}{c_N+\sqrt{c_N}z} + (c_n-1)\frac{s_{c_N}(z)}{\sqrt{c_N}}.\nonumber
% \end{align}
Here the contour $\oint_{\mathcal{C}}$ is closed and taken in the positive direction in the complex plane,  enclosing the support of $F^c(x)$. 
The main result is stated in the following theorem.
\begin{theorem}\label{cltmain}
    Under Assumptions \ref{asp1*}, \ref{asp2} and \ref{asp3}, 
%and after truncation of the
%data, 
let $f_1, f_2, \ldots, f_k$ be functions on $\mathbb{R}$ and analytic on an open interval containing $\Bigl[-2+\frac{1}{\sqrt{c_N}}, \; 2+\frac{1}{\sqrt{c_N}}\Bigr]$.
Then, the random vector $\left( G_{n}(f_1), \ldots, G_{n}(f_k)\right)$ forms a tight sequence in $n$ and converges weakly to a centered Gaussian vector $(X_{f_1}, \ldots, X_{f_k})$ with the covariance function
\begin{align*}
     Cov(X_f, X_g) &= -\frac{1}{4\pi^2} \oint_{\mathcal{C}_1}\oint_{\mathcal{C}_2} f(z_1)g(z_2)Cov(M(z_1),M(z_2))\dif z_1\dif z_1
\end{align*}
where 
\begin{align*}
    Cov(M(z_1),M(z_2))= 2 \left[\frac{s_c^{\prime}(z_1)s_c^{\prime}(z_2)}{\{s_c(z_1)-s_c(z_2)\}^2}-\frac{1}{(z_1-z_2)^2}\right] - \frac{2s_c^{\prime}(z_1)s_c^{\prime}(z_2)}{
\left\{1+s_c(z_1)/\sqrt{c}\right\}^2
\left\{1+s_c(z_2)/\sqrt{c}\right\}^2
    },
\end{align*}
and $s_c(z)$ is defined in \eqref{s_c(z)}.
\end{theorem}

\begin{remark}
    Theorem \ref{cltmain} establishes a unified CLT for LSS of $\mathbf{B}_n$ when $p / n \rightarrow c \in(0, \infty]$. This result is consistent with the results of $\bR_n$ when $p/n \rightarrow c \in(0, \infty)$ in Theorem 1 of \cite{gao2017high} and Theorem 3.2 of \cite{yin2023central}.
\end{remark}

    In particular, when $p/n\to\infty$, 
    $$
    G_n(f)=n \int_{-\infty}^{+\infty} f(x) \mathrm{d}\left\{F^{\mathbf{B}_n}(x)-F^{c_N}(x)\right\}+\frac{1}{2 \pi i} \oint_{\mathcal{C}}f(z) 
    \left\{
    \frac{s^3(z)+s(z)-s^{\prime}(z)s(z)}{s^2(z)-1}-\frac{1}{z}
    \right\}
    \mathrm{d} z,
    $$
    where $s(z)=\frac{-z+\sqrt{z^2-4}}{2}$ is the Stieltjes transform of the semicircle law. 
    Then we have the following result.

\begin{theorem}\label{ult-cltmain}
  With the same notations and assumptions given in Theorem \ref{cltmain} with $p \asymp n^t, t>1$, then the random vector $\left( G_{n}(f_1), \ldots, G_{n}(f_k)\right)$ forms a tight sequence in $n$ and converges weakly to a centered Gaussian vector $(X_{f_1}, \ldots, X_{f_k})$ with the covariance function
\begin{align*}
     Cov(X_f, X_g) &= -\frac{1}{4\pi^2} \oint_{\mathcal{C}_1}\oint_{\mathcal{C}_2} f(z_1)g(z_2)Cov(M(z_1),M(z_2))\dif z_1\dif z_1,
\end{align*}
where 
\begin{align}\label{covreal}
    Cov(M(z_1),M(z_2))= 2 \left[\frac{s^{\prime}(z_1)s^{\prime}(z_2)}{\{s(z_1)-s(z_2)\}^2}-\frac{1}{(z_1-z_2)^2}\right] - 2s^{\prime}(z_1)s^{\prime}(z_2).
\end{align}
\end{theorem}
\begin{remark}
   Theorem \ref{ult-cltmain} establishes a novel CLT for LSS of the renormalized sample correlation matrix $\mathbf{B}_n$ in the ultrahigh-dimensional regime where $p/n \to \infty$, which constitutes the main contribution of this paper. The proof technique is different from the classical case where $p/n\to c\in(0,\infty)$.
     In particular, we develop 
     a novel determinant equivalent form
for the resolvent of the renormalized correlation matrix, under ultrahigh dimensional context (see proof of Lemma \ref{le-J}). Theorem \ref{cltmain} provides a unified formulation of the limiting results for both $p/n \to c \in (0,\infty)$ and $p/n \to \infty$.
\end{remark}
\begin{Corollary}
With the same notations and assumptions given in Theorem \ref{cltmain}, consider three analytic functions $f_2(x)=x^2$,$f_3(x)=x^3$, $f_4(x)=x^4$, we have
$$
\begin{aligned}
& G_n\left(f_2\right)=\operatorname{tr}\left(\mathbf{B}_n^2\right)-n+2 \xrightarrow{d} \mathcal{N}(0,4),\\
&G_n\left(f_3\right)=\operatorname{tr}\left(\mathbf{B}_n^3\right)-\frac{n-4}{\sqrt{c_N}} \xrightarrow{d} \mathcal{N}\left(0, 6+\frac{36}{c}\right),\\
& G_n\left(f_4\right)=\operatorname{tr}\left(\mathbf{B}_n^4\right)-\left(\frac{(n-1)^2}{p}+2n-5-\frac{6}{c_N}\right)\xrightarrow{d} \mathcal{N}\left(0, 72+\frac{288}{c}+\frac{144}{c^2}\right).
\end{aligned}
$$
\end{Corollary}

\iffalse
 \begin{align}\label{covreal}
    Cov(M(z_1),M(z_2))= 2 [\frac{s^{'}(z_1)s^{'}(z_2)}{\{s(z_1)-s(z_2)\}^2}-\frac{1}{(z_1-z_2)^2}] - (|\Psi|^2+1)s^{'}(z_1)s^{'}(z_2),
 \end{align}
 under the real random-variable
 case, $\Psi=\frac{\Expe \{y_{11}-\Expe y_{11}\}^2}{\Expe |y_{11}-\Expe y_{11}|^2}=1$,
 $s(z)$ is the Stieltjes transform of the semicircle law. 
 When ${y_{ij}}$ are complex variables,
 $Cov(M(z_1),M(z_2))$ is given by
  assuming that $\Psi$ are the same for
 $i=1,2, \ldots, p$, 
 \begin{align}
     Cov(X_f, X_g) &= -\frac{1}{4\pi^2} \oint_{\mathcal{C}_1}\oint_{\mathcal{C}_2} f(z_1)g(z_2)Cov(M(z_1),M(z_2))\dif z_1\dif z_1
 \end{align}
 where
 \begin{align}
    Cov(M(z_1),M(z_2))&= \eqref{covreal}-[\frac{s^{'}(z_1)s^{'}(z_2)}{\{s(z_1)-s(z_2)\}^2}-\frac{1}{(z_1-z_2)^2}]\nonumber\\
     &+|\Psi|^2\frac{s^2(z_1)s^2(z_2)}{(-1+s^2(z_1))(-1+s^2(z_2))(-1+|\Psi|^2s(z_1)s(z_2))^2}.
 \end{align}
 \fi

\section{ Application of CLTs to hypothesis test}\label{ap}
In this section, we provide a statistical application of LSS for renormalized sample correlation matrix $\bB_n$. It is the independence test for high and ultra-high dimensional random vectors, namely the hypothesis
\[
    H_0:\; 
    \mathcal{R}=\bI_p\quad
    \text{vs.} \quad H_1:\;  \mathcal{R}\neq\bI_p.
\]
We aim to propose 
a test statistic that can work when  $p/n\to c\in (0,\infty]$.

Motivated by the Frobenius norm of $\mathcal{R}-\bI_p$ used in \cite{Schott}, \cite{gao2017high}
and \cite{yin2023central}, with the relationship
$$
\tr(\bR_n-\bI_p)^2=\frac{p}{N}(\tr \bB_n^2+p)-p,
$$
we consider the following test statistic constructed from the renormalized correlation matrix $\bB_n$,
\begin{align*}
    T := \tr\bB_n^2.
\end{align*}
We reject $H_0$ when $T$ is too large. 
%The test statistic is linked to particular forms of LSS of $\bB_{n}$ by taking $f(x)=x^2$, that is, 
%\[
%    T = \tr(\bB_{n}^2).
%\]
By 
by taking $f(x)=x^2$ in 
 Theorem \ref{cltmain}, the limiting null distribution of $T$ is given in the following theorem. 

\begin{theorem}\label{thm:limiting_null_dist}
  Suppose that Assumptions \ref{asp1*}, \ref{asp2} and \ref{asp3} hold, under $H_0$, we have
as $n\to\infty$,
$$
\frac{1}{2}(T-n+2)\xrightarrow{D}
\calN( 0,1 ).
$$
\iffalse
\begin{align*}
        \frac{1}{2}( T - \mu_T ) \to \calN( 0, \; 1 ),
    \end{align*}
    where 
    \begin{align*}
        \mu_T & = n\eta_2-\frac{1}{2 \pi i}\oint_{\mathcal{C}} z^2\Theta_n (s_{c_N}(z))\dif z,\\
        \eta_2 &=c_N\sum_{i=0}^2(-1)^i C_2^i(c_N)^{-1+i}\beta_{2-s}+(1-c_N)c_N,\\
       \beta_k&=\sum_{r=0}^{k-1} \frac{1}{r+1}\binom{k}{r}\binom{k-1}{r} c_N^r, k\geq 1, \ \ \beta_0=1.
    \end{align*}   
    \fi
\end{theorem}

    Theorem \ref{thm:limiting_null_dist} 
    establishes the unified CLT for $T$ under $H_0$ when $p / n \rightarrow c \in(0, \infty]$. Based on these, we employ the following procedure for testing the null hypothesis:
     \[
        \text{Reject $H_0$ if } \frac{1}{2}(T-n+2) > z_{\alpha},
    \]
    where $z_{\alpha}$ is the upper-$\alpha$ quantile of the standard normal distribution at nominal level $\alpha$.

\section{Simulations}

In this section, we implement some simulation studies to examine
\begin{itemize}
     \item[(1)] the LSD of the renormalized sample correlation matrix $\bB_{n}$;%, as stated in Theorem \ref{thm:LSD};
    \item[(2)] finite-sample properties of some LSS for $\bB_n$ by comparing their empirical means and
variances with theoretical limiting values;%, as stated in Theorem \ref{cltmain};

\item[(3)] finite-sample performance of independence test.
\end{itemize}

\subsection{Limiting spectral distribution}

In this section, simulation experiments are conducted to verify the LSD of the renormalized sample correlation matrix $\bB_{n}$, as stated in Theorem \ref{thm:LSD}. We generate data $y_{ij}$ from three populations, drawing histograms of eigenvalues of $\bB_{n}$ and comparing them with theoretical densities. Specifically, three types of distributions for $y_{ij}$ are considered:
\begin{itemize}
% \begin{enumerate}
	\item[(1)] $y_{ij}$ follows the standard normal distribution;
	\item[(2)] $y_{ij}$ follows the exponential distribution with rate parameter $2$;
  \item[(3)] $y_{ij}$ follows the Poisson distribution with parameter $1$.
% \end{enumerate}
\end{itemize}

\begin{figure}[htbp]
	\centering
\begin{subfigure}{.32\textwidth}
		\centering		\includegraphics[width=.9\linewidth]{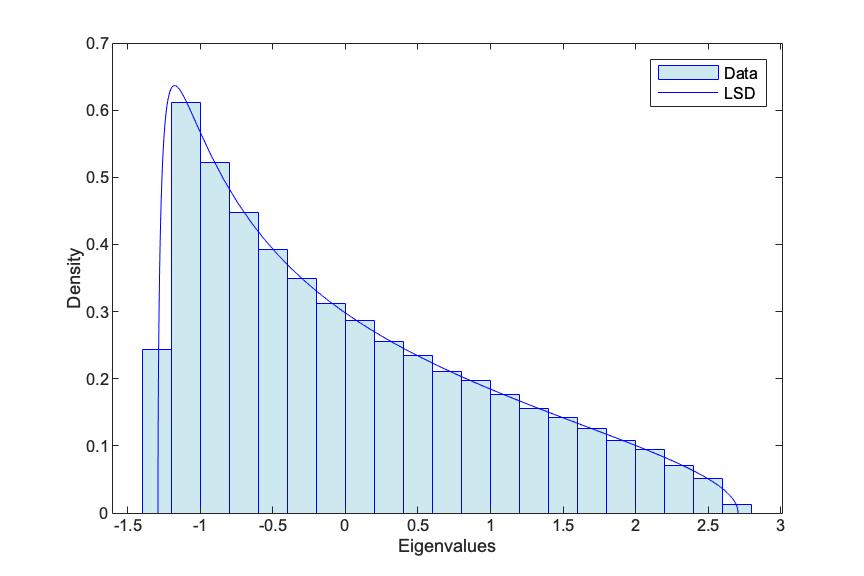}
		\caption{N(0,1)}
	\end{subfigure}%
	\begin{subfigure}{.32\textwidth}
		\centering		\includegraphics[width=.9\linewidth]{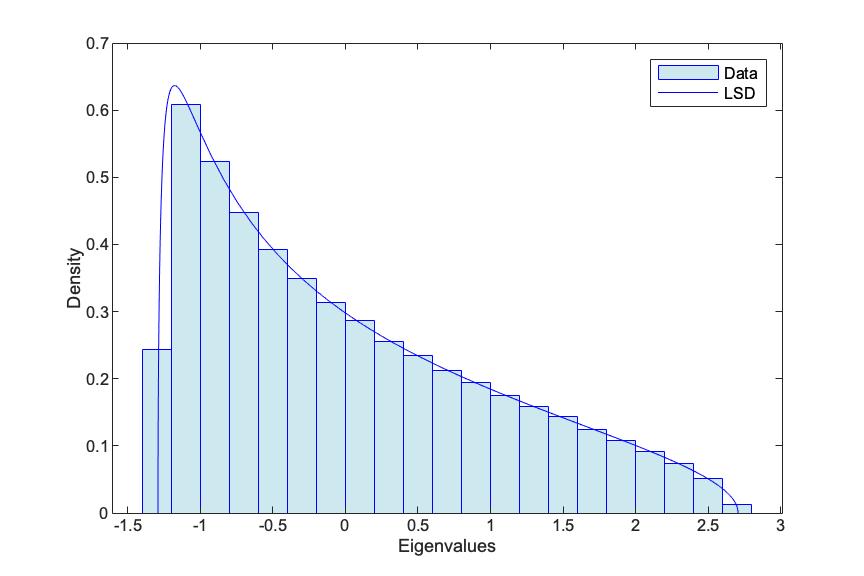}
		\caption{Exponential(2)}
	\end{subfigure}
	\begin{subfigure}{.32\textwidth}
		\centering		\includegraphics[width=.9\linewidth]{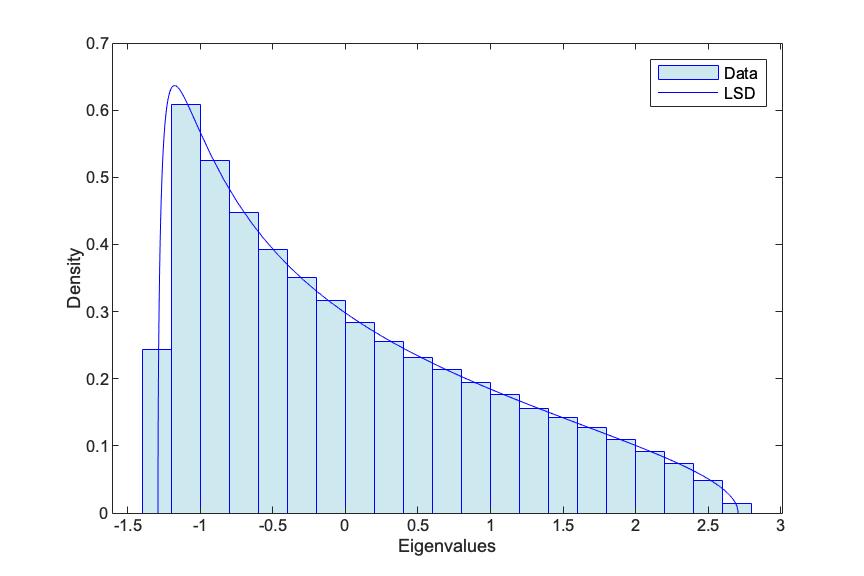}
		\caption{Poisson(1)}
	\end{subfigure}%
 \\
	\begin{subfigure}{.32\textwidth}
		\centering		\includegraphics[width=.9\linewidth]{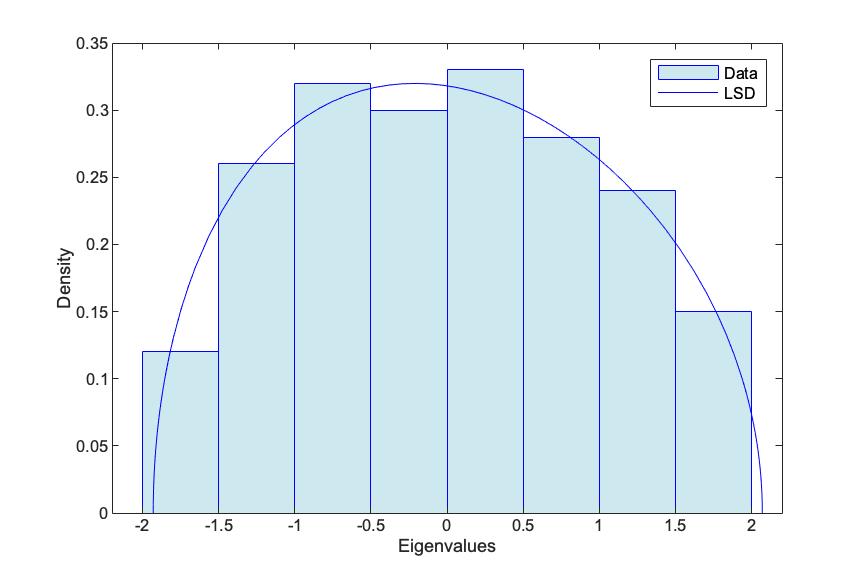}
		\caption{N(0,1)}
	\end{subfigure}%
	\begin{subfigure}{.32\textwidth}
		\centering		\includegraphics[width=.9\linewidth]{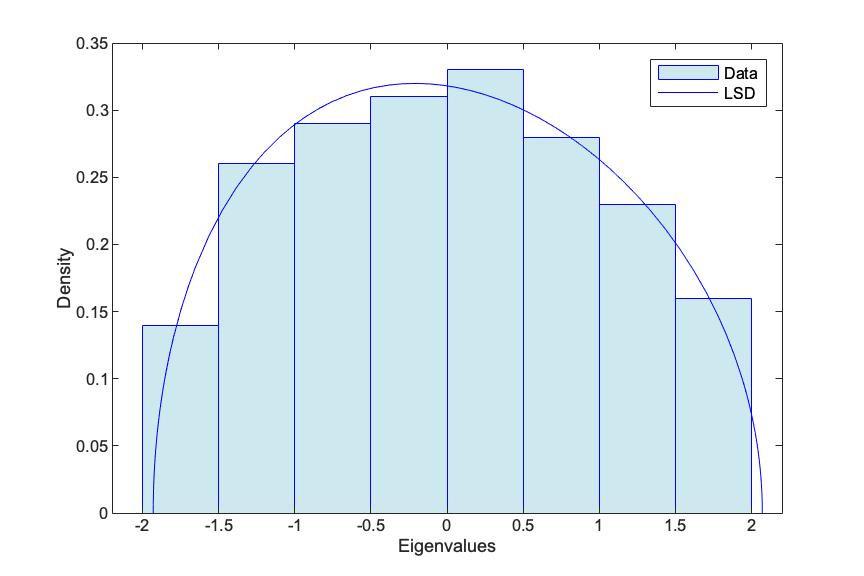}
		\caption{Exponential(2)}
	\end{subfigure}
	\begin{subfigure}{.32\textwidth}
		\centering		\includegraphics[width=.9\linewidth]{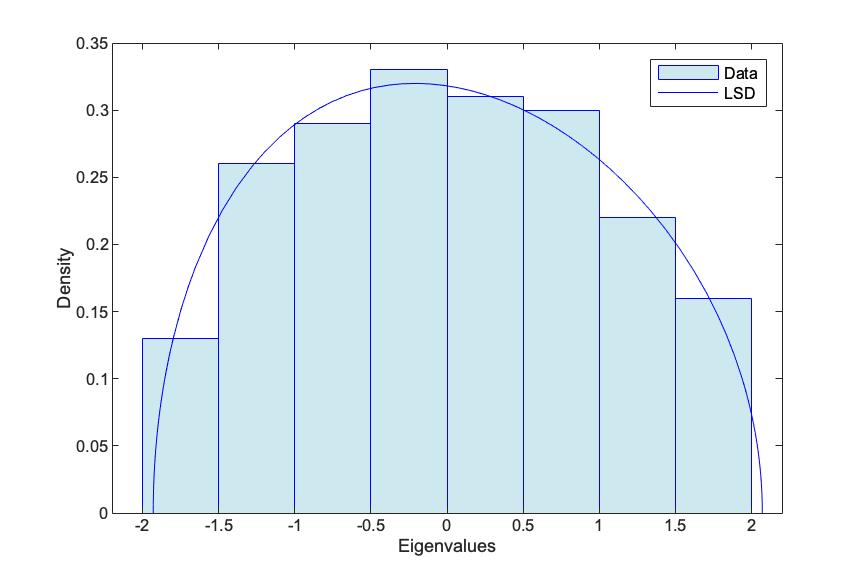}
		\caption{Poisson(1)}
	\end{subfigure}%
 \\
 	\begin{subfigure}{.32\textwidth}
		\centering		\includegraphics[width=.9\linewidth]{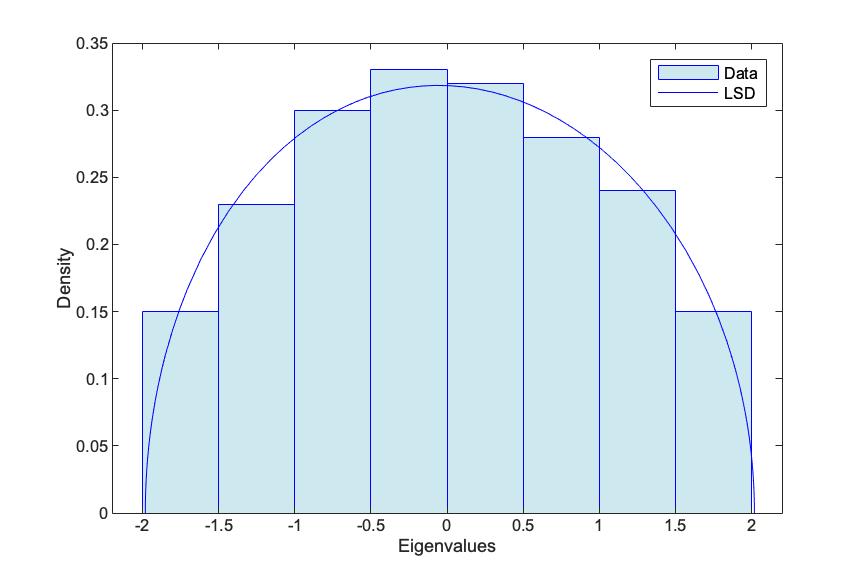}
		\caption{N(0,1)}
	\end{subfigure}%
	\begin{subfigure}{.32\textwidth}
		\centering
        \includegraphics[width=.9\linewidth]{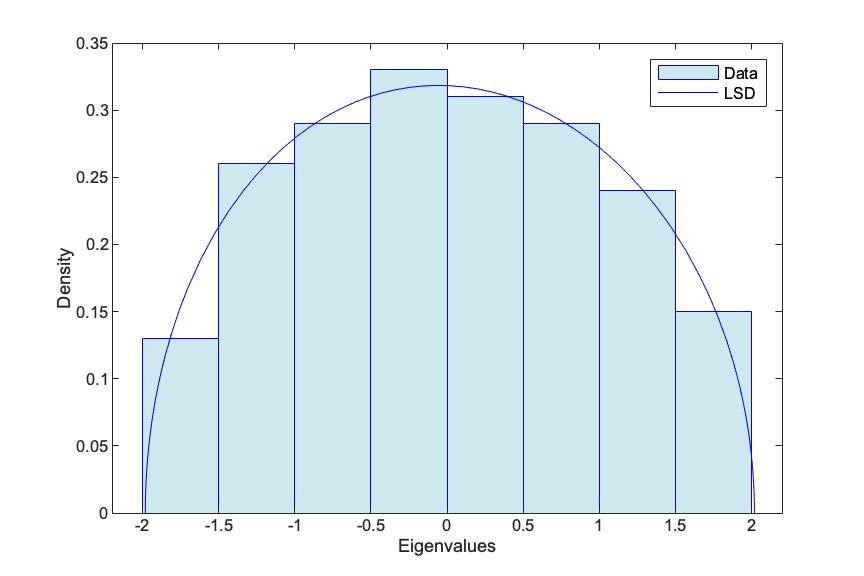}
		\caption{Exponential(2)}
	\end{subfigure}
	\begin{subfigure}{.32\textwidth}
		\centering
		\includegraphics[width=.9\linewidth]{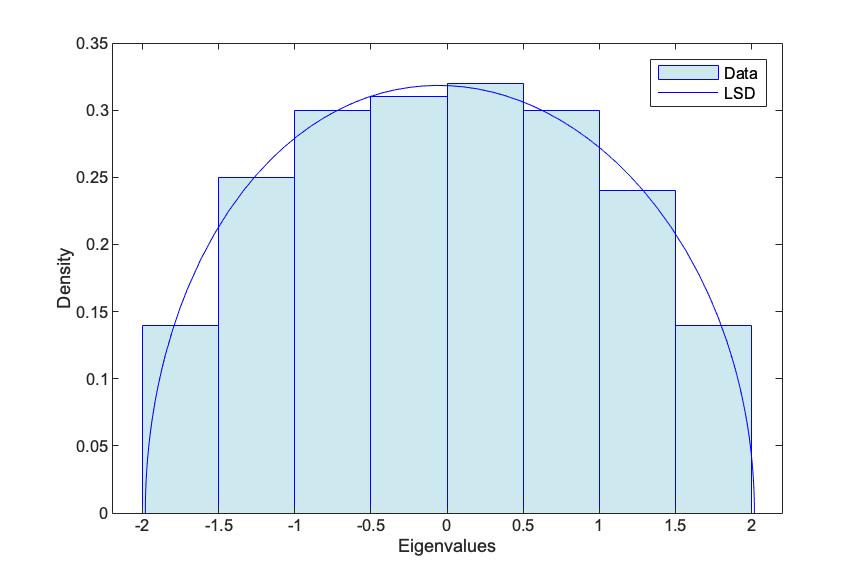}
		\caption{Poisson(1)}
	\end{subfigure}%
	\caption{Histograms of sample eigenvalues of $\bB_{n}$, fitted by  LSD (blue solid curves).
 In the first row, $(p,n)=(10^4,5000)$, 
 in the second row,
 $(p,n)=(200^2,200)$, in the third row $(p,n)=(200^{2.5},200)$.
 }\label{lsd1}
\end{figure}

The dimensional settings are $(p, n)=(10^4, 5000),(200^2,200),(200^{2.5},200).$
We display histograms of eigenvalues of $\bB_{n}$ generated by three populations under various $(p, n)$ in Figure \ref{lsd1}. This reveals that all histograms align with their LSD, affirming the accuracy of our theoretical results.

%\iffalse
%\begin{figure}[htbp]
%	\centering
%	\begin{subfigure}{.32\textwidth}
%		\centering		\includegraphics[width=.9\linewidth]{ult-corr/figure/LSD-Gaussian-n=200-p=n25.jpg}
%		\caption{N(0,1)}
%	\end{subfigure}%
%	\begin{subfigure}{.32\textwidth}
%		\centering
%        \includegraphics[width=.9\linewidth]{ult-corr/figure/LSD-exp-n=200-p=n25.jpg}
%		\caption{Exponential(2)}
%	\end{subfigure}
%	\begin{subfigure}{.32\textwidth}
%		\centering
%		\includegraphics[width=.9\linewidth]{ult-corr/figure/LSD-po-n=200-p=n25.jpg}
%		\caption{Poisson(1)}
%	\end{subfigure}%
%	\caption{Histograms of sample eigenvalues of $\bB_{n}$ with $(p,n)=(200^{2.5},200)$. The curves are density functions of their corresponding limiting spectral distribution.}\label{lsd2}
%\end{figure}
%\fi

\subsection{CLT for LSS}

In this section, we implement some simulation studies to examine finite-sample properties of some LSS for $\bB_{n}$ by comparing their empirical means and variances with theoretical limiting values, as stated in Theorem \ref{cltmain}. In the following, we present the numerical simulation of CLT for LSS. Let $\bar{G}_{n}(f_r)=\frac{G_{n}(f_r)}{\sqrt{\Var(X_{f_r})}}$. First, we examine $\bar{G}_{n}(f_r)\stackrel{d}{\to}N(0,1)$, $f_r=x^r(r = 2,3,4)$, by Theorem \ref{cltmain}. Two types of data distribution of $y_{ij}$ are consider:
\begin{itemize}
% \begin{enumerate}
	%\item $y_{ij}$ follows the standard normal distribution;
	%\item $y_{ij}$ follows the Chi-squared distribution with degree of freedom $1$. 
 \item[(1)] Gaussian data: $y_{ij}\sim N(0,1)$ i.i.d. for $1\leq i\leq p$, $1\leq j \leq n$;
 \item[(2)] Non-Gaussian data: $y_{ij}\sim \chi^2(2)$ i.i.d. for $1\leq i\leq p$, $1\leq j \leq n$.
% \end{enumerate}
\end{itemize}

Empirical mean and variance of $\{\bar{G}_{n}(f_r)\}$, $f_r=x^r, r=2,3,4$, are calculated for various  $c_n$.
The dimensional settings are $(p,n)=(1000,500),(300^2,300),(500^2,500),(100^{2.5},100)$ with $c_n=2,300,500,1000$.
%For each pair of $(p,n)$, $5000$ independent replications are used to obtain the empirical values. Tables \ref{cltlssp2n}-\ref{cltlsspn2.5} report the empirical values of $\{\bar{G}_{n}(f_r)\}$ with $p=2n$, $p=n^2$ and  $p=n^{2.5}$ respectively. 
As shown in Tables \ref{clt}, the empirical mean and variance of $\{\bar{G}_{n}(f_r)\}$ perfectly match their theoretical limits 0 and 1 under all scenarios.

\begin{table}[htbp]
\caption{Empirical mean and variance of $\bar{G}_n\left(f_i\right), i=2,3,4$ from 5000 replications with $c_n=2,300,500,1000$. Theoretical mean and variance are 0 and 1, respectively.}
\centering
%\renewcommand{\arraystretch}{1.3}
\iffalse
\begin{tabular}{ccccccccc}
\toprule
 & \multicolumn{2}{c}{$\bar{G}_n\left(f_2\right)$} 
 &    
 & \multicolumn{2}{c}{$\bar{G}_n\left(f_3\right)$} 
 &
 & \multicolumn{2}{c}{$\bar{G}_n\left(f_4\right)$} \\
\cmidrule{2-3} \cmidrule{5-6} \cmidrule{8-9} 
 $c_n$    & mean & var& & mean & var & &mean & var \\
\midrule
          &       &       &       & \multicolumn{5}{l}{\textit{Gaussian data}} \\
2   & 0.0180 & 1.0079 & &  -0.0504 & 0.9737& & 0.0642 & 0.9793 \\
300  & 0.0369 & 0.9974 & &-0.2273 & 0.9973 & &0.0320 & 0.9914 \\
500  & 0.0285 & 0.9837 & &-0.2022 & 1.0033 & &-0.0649 & 0.9718 \\
1000  & 0.0289 & 0.9889 & &-0.1143 & 0.9770 & &-0.0300 & 0.9936 \\
& & & & \textit{Non-Gaussian  data} \\
2   &  0.0567 & 1.1201 && -0.0596 & 1.0672 & &0.0075 & 1.0342\\
300  & -0.0714 &1.0326 &&-0.2417 & 1.0599 && 0.0383&1.0110 \\
500  & -0.0131& 1.0240 && -0.1709 & 1.0236 && -0.0055 & 1.0260 \\
1000  &0.0040 & 1.0840 && -0.1431 & 0.9842 && 0.0219 & 1.0474 \\
\bottomrule
\end{tabular}
\fi
\begin{tabular}{ccccccccc}
\toprule
 & \multicolumn{2}{c}{$\bar{G}_n\left(f_2\right)$} 
 &    
 & \multicolumn{2}{c}{$\bar{G}_n\left(f_3\right)$} 
 &
 & \multicolumn{2}{c}{$\bar{G}_n\left(f_4\right)$} \\
\cmidrule{2-3} \cmidrule{5-6} \cmidrule{8-9} 
 $c_n$    & mean & var& & mean & var & &mean & var \\
\midrule
          &       &       &       & \multicolumn{5}{l}{\textit{Gaussian data}} \\
2   & 0.0090 & 1.0079 & &  -0.0103 & 0.9737& & 0.0040 & 0.9793 \\
300  & 0.0185 & 0.9974 & & -0.0919 & 0.9777 & & 0.0037 & 0.9785 \\
500  & 0.0143 & 0.9837 & & -0.0821 & 0.9914 & & -0.0076 & 0.9639 \\
1000  & 0.0144 & 0.9889 & & -0.0465 & 0.9712 & & -0.0035 & 0.9896 \\
& & & & \textit{Non-Gaussian  data} \\
2   & 0.0284 & 1.1201 && -0.0122 & 1.0672 & & 0.0005 & 1.0342\\
300  & -0.0357 &1.0326 && -0.0977 & 1.0390 && 0.0045 & 0.9974 \\
500  & -0.0066 & 1.0240 && -0.0694 & 1.0112 && -0.0006 & 1.0179 \\
1000  & 0.0020 & 1.0840 && -0.0582 & 0.9785 && 0.0026 & 1.0432 \\
\bottomrule
\end{tabular}
\label{clt}
\end{table}

\subsection{Hypothesis test}

Numerical simulations are conducted to find empirical size and powers of our proposed test statistic. The random variables $(x_{ij})$ are generated from:
\iffalse
    To evaluate the finite sample performance of these test statistics, data are generated from various model scenarios. For each scenario, We run $5000$ independent replications. 
    To examine empirical size, we consider two scenarios with respect to the data $\bY=(y_{ij})_{p\times n}$:
    \fi
    \begin{itemize}
        \item[(1)] Gaussian data: $x_{ij}\sim N(0,1)$ i.i.d. for $1\leq i\leq p$, $1\leq j \leq n$;
        \item[(2)] Non-Gaussian data: $x_{ij}\sim (\chi^2(2)-2)/2$ i.i.d. for $1\leq i\leq p$, $1\leq j \leq n$.
    \end{itemize}
    And we consider the following two settings of $\bSigma$:
 \iffalse
 To examine their empirical power, we consider the following two alternatives:
 A matrix $\bX$ with i.i.d. components is generated first following Gaussian/chi-squared models, then the data matrix $\bY=\bSigma\bX$
    is constructed with a $p\times p$ covariance matrix $\bSigma$ defined as follows:
    \fi
    \begin{itemize}

        \item $\bSigma_1=\left(s_{i, j, \theta}\right)_{p \times p}, s_{i, j, \theta}=\delta_{\{i=j\}}+\delta_{\{i \neq j\}} \theta^{|i-j|}$, $ i, j=1, \ldots, p$,
 \item $\bSigma_2=\left(s_{i, j, \eta}\right)_{p \times p}, s_{i, j, \eta}=\delta_{\{i=j\}}+\delta_{\{i \neq j\}} \eta, i, j=1, \ldots, p$,
    \end{itemize}
where $\theta, \eta$ are two parameters satisfying $|\theta|<1,0<\eta<1$. The parameter setting is as follows:

\begin{itemize}
\item $\theta=\eta=0$ to evaluate empirical size;
\item $\theta=0.20,0.25$  to evaluate empirical power of $\bSigma_1$;
    \item $\eta=0.007, 0.011$ to evaluate empirical power of $\bSigma_2$.
\end{itemize}

Table \ref{tab:size_power} reports the empirical size and power for different $c_n$. The dimensional settings are $(p, n)=(1200,600),(50^2, 50),(100^2, 100),(200^2, 200)$ with $c_n=2,50,100,200$, and the nominal significance level is fixed at $\alpha=0.05$. This shows our test statistic is robust in both high and ultra-high dimensional settings and performs stably for Gaussian and non-Gaussian data.

\begin{table}[htbp]
  \centering
  \caption{Empirical size and power from $5000$ replications for Gaussian and Non-Gaussian data with different $c_n$.}
  \begin{tabular}{ccccrccccccc}
    \toprule
          & Size  &     & \multicolumn{2}{c}{Power of $\bSigma_1$} &            & \multicolumn{2}{c}{Power of $\bSigma_2$} \\
    \cmidrule{2-2} \cmidrule{4-5} \cmidrule{7-8} 
    $c_n$   & $\theta=\eta=0$ &   & $\theta=0.20$ & $\theta=0.25$ 
    && $\eta=0.007$
    & $\eta=0.011$
    \\
    \midrule
          &       &       &       & \multicolumn{4}{l}{\textit{Gaussian data}} \\
    2 & 0.0528 &  & 0.9970 & 1
    && 0.5954 & 0.9866\\
    50 & 0.044 &  & 0.608 & 0.902 
    && 0.7688 & 0.9878\\
    100 & 0.0456 &   & 0.9884 & 1 
    &&0.9999 &1\\
    200 & 0.0512 &   & 1 & 1 
     &&1 &1\\
          &       &       &       & \multicolumn{4}{l}{\textit{Non-Gaussian  data}}
          \\
    2 & 0.0498 &  & 0.9964 & 1
    && 0.5908 & 0.9814\\
    50 & 0.06 &  & 0.6278 & 0.922 
     &&0.7372 &0.98\\
    100 & 0.0542 &   & 0.9878 & 1 
    &&0.9997&1\\
    200 & 0.055 &   & 1 & 1
      &&1 &1\\
    \bottomrule
  \end{tabular}
  \label{tab:size_power}
\end{table}

\section{Proofs}
\subsection{Notations}
The following notations are used throughout the proofs. Let
$\mathbf{Y}^{\top} = (\tilde{\by}_1,\ldots,\tilde{\by}_p),$
then $\bB_n$ can be written as
$$
\bB_n=\sqrt{\frac{p}{n-1}}\left(\frac{n-1}{p} \tilde{\bY}_n \tilde{\bY}_n^{\top}-\boldsymbol{\Phi}\right),~
\tilde{\mathbf{Y}}_n=\left(\frac{\bPhi\tilde{\mathbf{y}}_1}{\left\|\bPhi\tilde{\mathbf{y}}_1\right\|}, \frac{\bPhi\tilde{\mathbf{y}}_2}{\left\|\bPhi\tilde{\mathbf{y}}_2\right\|}, \ldots, \frac{\bPhi\tilde{\mathbf{y}}_p}{\left\|\bPhi\tilde{\mathbf{y}}_p\right\|}\right).
$$
Denote
\begin{align*}
\bA_n&=\sqrt{\frac{p}{N}}\left(\frac{N}{p} \underline\bR_n-\bI_n\right),\ \underline\bR_n=\tilde{\bY}_n \tilde{\bY}_n^{\top},\ N=n-1,\ 
c_n=p/n,\
c_N=p/N,\\
s_n(z)&=\frac{1}{n} \operatorname{tr}\left(\bA_n-z \bI_n\right)^{-1},\ s_n^{\bB_n}(z)=\frac{1}{n} \operatorname{tr}\left(\bB_n-z \bI_n\right)^{-1}, \quad z\in \mathbb{C}^{+},\\
% \underline{m}_n(z)=\frac{1}{n} \operatorname{tr}\left(\underline{\bR}_n-z \bI_n\right)^{-1},
\tilde{\mathbf{Y}}_n&=\left(\frac{\boldsymbol{\Phi} \tilde{\mathbf{y}}_1}{\left\|\boldsymbol{\Phi} \tilde{\mathbf{y}}_1\right\|}, \frac{\boldsymbol{\Phi} \tilde{\mathbf{y}}_2}{\left\|\boldsymbol{\Phi} \tilde{\mathbf{y}}_2\right\|}, \ldots, \frac{\boldsymbol{\Phi} \tilde{\mathbf{y}}_p}{\left\|\boldsymbol{\Phi} \tilde{\mathbf{y}}_p\right\|}\right)=\left(\br_1,\ldots,\br_p\right),\
 \tilde{\bY}_k=\left(\br_1, \cdots, \br_{k-1}, \br_{r+1}, \cdots, \br_p\right),\\
\underline\bR_{nk}&=\tilde{\bY}_k\tilde{\bY}_k^{\top},\ 
\bA_{nk}=\sqrt{\frac{p}{N}}\left(\frac{N}{P} \underline\bR_{nk}-\bI_n\right),\
\bA_{nkj}=\bA_{nk}-\sqrt{\frac{N}{p}}\br_j\br_j^{\top},\\
\bQ(z)&=\bA_n-z\bI_n,\
\bQ_k(z)=\bA_{nk}-z\bI_n,\
\bQ_{kj}(z)=\bA_{nkj}-z\bI_n,
\\
\beta_k(z)&=\frac{1}{\sqrt{c_N}+\br_k^{\top}\bQ_k^{-1}(z)\br_k},\
\tilde\beta_k(z)=\frac{1}{\sqrt{c_N}+\operatorname{tr}\bQ_k^{-1}(z)/n},\
b_n(z)=\frac{1}{\sqrt{c_N}+\E\operatorname{tr}\bQ_k^{-1}(z)/n},\\
b_1(z)&=\frac{1}{\sqrt{c_N}+\E\operatorname{tr}\bQ_{12}^{-1}(z)/n},\
\gamma_k(z)=\br_k^{\top} \bQ_k^{-1} (z)\br_k-\E\frac{1}{n} \operatorname{tr}  \bQ_k^{-1}(z),
\beta_{kj}(z)=\frac{1}{\sqrt{c_N}+\br_k^{\top}\bQ_{kj}^{-1}(z)\br_k},
\\
\varepsilon_k(z)&=\br_k^{\top} \bQ_k^{-1}(z) \br_k-\frac{1}{n} \operatorname{tr}  \bQ_k^{-1}(z),\
\delta_k(z)=\br_k^{\top} \bQ_k^{-2}(z) \br_k-\frac{1}{n} \operatorname{tr}  \bQ_k^{-2}(z).
\end{align*}
We denote by $K$ some constant which may take different values at different appearances.

By the results in \cite{bai2004clt}, we have $\|\bQ_k(z)^{-1}\|\leq K$, 
$\left|\operatorname{tr}\left(\bQ^{-1}(z)-\bQ_k^{-1}(z)\right) \bM\right| \leq \|\bM\|c_n^{-\frac{1}{2}}$, 
$|\beta_k(z)|\leq Kc_n^{-\frac{1}{2}}$, $|\tilde\beta_k(z)|\leq Kc_n^{-\frac{1}{2}}$, $|b_n(z)|\leq Kc_n^{-\frac{1}{2}}$. And straightforward calculation gives  
\begin{align}\label{expand-beta}
&\bQ^{-1}(z)-\bQ^{-1}_k(z)=-\bQ_k^{-1}(z)\br_k\br_k^{\top}\bQ_k^{-1}(z)\beta_k(z), \nonumber\\
&\beta_k(z)=b_n(z)-b_n(z)\gamma_k(z)\beta_k(z)=\tilde\beta_k(z)-\tilde\beta_k(z)\epsilon_k(z)\beta_k(z).
\end{align}

\subsection{Proof of Theorem \ref{ult-thm:LSD}}\label{se:proof-lsd}
Since
\begin{align}\label{sb-sn}
s_n^{\bB_n}(z)=s_n(z)-\frac{1}{n}\frac{\sqrt{c_N}}{z(\sqrt{c_N}+z)},
\end{align}
for all $z \in \mathbb{C}^{+}$, the difference between $s_n^{\mathbf{B}_n}(z)$ and $s_n(z)$ is a deterministic term of order $O(1 / n)$. Therefore, to show that $s_n^{\mathbf{B}_n}(z) \rightarrow s(z)$ almost surely, it suffices to prove that $s_n(z) \rightarrow s(z)$ almost surely. Here $s(z)$ is the Stieltjes transform of semicircular law \eqref{eq:sc_density}. We now proceed with the proof in the following four steps:
\iffalse
Since 
\begin{align}\label{sb-sn}
s_n^{\bB_n}(z)=s_n(z)-\frac{1}{n}\frac{\sqrt{c_N}}{z(\sqrt{c_N}+z)},
\end{align}
the LSD of $\bB_n$ can be obtained directly by \eqref{mp-lsd} when $p/n\to c\in(0,\infty)$. We thus focus on the case where $p/n\to \infty$.
By using the relationship \eqref{sb-sn} again,
for all $z \in \mathbb{C}^{+}$, to prove $s_n^{\bB_n}(z) \rightarrow s(z)$ almost surely,   it suffices to prove  $s_n(z) \rightarrow s(z)$ almost surely.
Here $s(z)$ is the Stieltjes transform of semicircular law \eqref{eq:sc_density}.
We shall then proceed in our proof by taking the following four steps:
\fi
\begin{itemize}
\item[\textbf{Step 1.}] Truncation, centralization, and rescaling.
    \item[\textbf{Step 2.}] For any fixed $z \in \mathbb{C}^{+}=\{z, \Im(z)>0\}, s_n(z)-\mathbb{E} s_n(z) \rightarrow 0$, a.s..
    
\item[\textbf{Step 3.}] For any fixed $z \in \mathbb{C}^{+}, \mathbb{E}s_n(z) \rightarrow s(z)$.

\item[\textbf{Step 4.}] Outside a null set, $s_n(z) \rightarrow s(z)$ for every $z \in \mathbb{C}^{+}$.
\end{itemize}

Then, it follows that, except for this null set, $F^{\bB_n} \rightarrow F$ weakly, where $F$ is the distribution function of semicircular law in \eqref{eq:sc_density}.

\textbf{Step 1.}
Truncation, centralization, and rescaling.
By the moment condition $\E|\bx_{11}|^4<\infty$,  one may choose a positive sequence of $\{\Delta_n\}$ such that
$$\Delta_n^{-4} \mathbb{E}\left|x_{11}\right|^4 I_{\left\{\left|x_{11}\right| \geqslant \Delta_n \sqrt[4]{n p}\right\}} \rightarrow 0, \quad \Delta_n  \rightarrow  0, \quad \Delta_n \sqrt[4]{n p} \rightarrow \infty.$$
Recall $\mathbf{X}=\left(\mathbf{x}_1, \ldots, \mathbf{x}_n\right)_{p \times n}=\left(x_{i j}\right)$. Then we can write $\bB_n=\bB_n(x_{ij})=
\bPhi\bB_{n0}\bPhi$, where $$
\bB_{n0}=\sqrt{\frac{p}{N}}\left[\frac{1}{p}\bX^{\top}\bD_n\bX-\bI_n\right],~
\bD_n=\operatorname{Diag}\left(\frac{1}{s_{11}},\frac{1}{s_{22}},\ldots,\frac{1}{s_{pp}}\right),~
s_{kk}=\frac{1}{N}\be_k^{\top}\bX\bPhi\bX^{\top}\be_k,~k=1,\ldots,p.
$$
Let $\hat\bB_n=\hat\bB_n(\hat x_{ij})$, $\check\bB_n=\check\bB_n(\check x_{ij})$ and 
$\tilde\bB_n=\tilde\bB_n(\tilde x_{ij})$
be defined similarly to $\bB_n$ with $x_{ij}$ replaced by $\hat{x}_{ij}$, $\check{x}_{ij}$ and 
$\tilde{x}_{ij}$ respectively,  where $\hat x_{ij}=x_{ij}I_{\{|x_{ij}|\leq \Delta_n \sqrt[4]{n p}\}}$, $\check x_{ij}=\hat x_{ij}-\E\hat x_{ij}$, 
and $\tilde{x}_{ij}=\check x_{ij}/\sigma_n$ with $\sigma_n^2=\E|\hat x_{ij}-\E\hat x_{ij}|^2$.
And similarly define $\hat\bD_n$, $\check\bD_n$, $\tilde\bD_n$ and $\hat\bB_{n0}$, $\check\bB_{n0}$, 
$\tilde\bB_{n0}$. Note that $\hat\bD_n=\check\bD_n$ and $\check\bB_{n0}=\tilde\bB_{n0}$. Then by Theorems A.43-A.44 in \cite{BSbook} and Bernstein's inequality, we have
\begin{align*}
\|F^{\bB_n}-F^{\bB_{n0}}\|&\leq \frac{1}{n}\operatorname{rank}(\bB_n-\bB_{n0})\leq \frac{K}{n},\\
\|F^{\bB_{n0}}-F^{\hat\bB_{n0}}\|&\leq \frac{1}{n}\operatorname{rank}\left(\bX^{\top}\bD_n^{\frac{1}{2}}
-\hat\bX^{\top}\hat\bD_n^{\frac{1}{2}}\right)\leq
\frac{1}{n}\sum_{i=1}^p\sum_{j=1}^{n}I_{\{|x_{ij}|\geq \Delta_n \sqrt[4]{n p}\}}\to 0 \quad a.s.,\\
\|F^{\hat\bB_{n0}}-F^{\tilde\bB_{n0}}\|&=
\|F^{\hat\bB_{n0}}-F^{\check\bB_{n0}}\|
\leq 
\frac{1}{n}\operatorname{rank}\left(
\hat\bX^{\top}\hat\bD_n^{\frac{1}{2}}
-\check\bX^{\top}\check\bD_n^{\frac{1}{2}}
\right)=\frac{1}{n}\operatorname{rank}
(\E \hat\bX^{\top}\hat\bD_n^{\frac{1}{2}})=\frac{1}{n}.
\end{align*}
Thus in the rest of the proof of Theorem \ref{thm:LSD},  we assume $$\left|x_{i j}\right| \leqslant \Delta_n \sqrt[4]{n p}, \quad \mathbb{E} x_{i j}=0, \quad \mathbb{E}\left|x_{i j}\right|^2=1, \quad \mathbb{E}\left|x_{i j}\right|^4=\kappa+o(1)<\infty .$$

\textbf{Step 2.} Almost sure convergence of the random part. Let $\Expe_0(\cdot)$ denote expectation and $\Expe_j(\cdot)$ denote conditional expectation with respect to the $\sigma$-field generated by $\br_1, \br_2, \ldots, \br_p$, where $j=1,2,\ldots,p$. 
By Lemma 2.7 in \cite{Bai1998No} and Lemma 5 in \cite{gao2017high}, we can obtain for $q>2$,
\begin{align}\label{estimate-lsd}
&\E\left|\varepsilon_k(z)\right|^q\leq K\left(n^{-q/2}+n^{-q/2} p^{q/2-1} \Delta_n^{2q-4}\right),\
\E\left|\delta_k(z)\right|^q\leq K\left(n^{-q/2}+n^{-q/2} p^{q/2-1} \Delta_n^{2q-4}\right),
\nonumber\\
&\E\left|\tilde\beta_k(z)-b_n(z)\right|^q=O(n^{q/2}p^{-q}),\
\left|b_n(z)-b_1(z)\right|=O(n^{1/2}p^{-2/3}),\
\E\left|b_n(z)-\E\beta_k(z)\right|=O(np^{-2}),\nonumber\\
&\E\left|\gamma_k(z)-\varepsilon_k(z)\right|^q=O(n^{-q/2}),\ \E\left|\frac{1}{n} \operatorname{tr}\left(\bQ^{-1}(z)\bM\right)-\E \frac{1}{n} \operatorname{tr}\left(\bQ^{-1}(z) \bM\right)\right|^q=O(-n^{q/2}).
\end{align}

Write
\begin{align*}
s_n(z)-\mathbb{E} s_n(z)=-\frac{1}{n}\sum_{j=1}^p\left(\mathbb{E}_j-\mathbb{E}_{j-1}\right) \beta_j(z) \mathbf{r}_j^{\top}\mathbf{Q}_j^{-2}(z) \mathbf{r}_j.
\end{align*}
By using Lemma 2.1 in \cite{bai2004clt}, we have \begin{align*}
\E|s_n(z)-\mathbb{E} s_n(z)|^4\leq \frac{K}{n^4}\E\left(\sum_{j=1}^p\left|\left(\mathbb{E}_j-\mathbb{E}_{j-1}\right)\beta_j(z) \mathbf{r}_j^{\top} \mathbf{Q}_j^{-2}(z) \mathbf{r}_j\right|^2\right)^2
=O(n^{-2}),
\end{align*}
where in the last step, we use the fact that 
$\left|\beta_j(z)\right| \leq K c_n^{-\frac{1}{2}}$
and $\E\left|
\mathbf{r}_k^{\top} \mathbf{Q}_k^{-2}(z) \mathbf{r}_k
\right|^2\leq \E\left|
\delta_k(z)
\right|^2+\E \left|\frac{1}{n} \operatorname{tr} \mathbf{Q}_k^{-2}(z)\right|^2=O(1)$ by \eqref{estimate-lsd}. Therefore, we obtain $s_n(z)-\E s_n(z)=o_{a.s.}(1)$.

\textbf{Step 3.}  Convergence of the expected Stieltjes transform. 
Similarly to the proof of Lemma \ref{step3} in the Supplementary Material, and by applying the estimates in \eqref{estimate-lsd}, we obtain $$n\left[\E s_n(z)-s_{c_N}(z)\right] =O(1),$$ which implies that  $\E s_n(z)=s_{c_N}(z) +O(n^{-1})$. The details are omitted here. Moreover, since $s_{c_N}(z)=s(z)+o(1)$, it follows that $\E s_n(z)=s(z)+o(1)$.

\textbf{Step 4.} Completion of the proof of Theorem \ref{thm:LSD}. By Steps 2 and 3, for any fixed $z \in \mathbb{C}^{+}$, we have
$
s_n(z) \rightarrow s(z) \text {, a.s.. }
$
That is, for each $z \in \mathbb{C}^{+}$, there exists a null set $N_z$ (i.e., $P\left(N_z\right)=0$ ) such that
$
s_n(z, \omega) \rightarrow s(z) \text { for all } \omega \in N_z^c.$
Now, let $\mathbb{C}_0^{+}=\left\{z_m\right\}$ be a dense subset of $\mathbb{C}^{+}$ (e.g., all $z$ of rational real and imaginary parts) and let $N=\cup N_{z_m}$. Then
$$
s_n(z, \omega) \rightarrow s(z) \text { for all } \omega \in N^c \text { and } z \in \mathbb{C}_0^{+}.
$$
Let $\mathbb{C}_m^{+}=\left\{z \in \mathbb{C}^{+}, \Im z>1 / m,|z| \leq m\right\}$. When $z \in \mathbb{C}_m^{+}$, we have $\left|s_n(z)\right| \leq m$ . By Vitali’s convergence theorem, we have
$$
s_n(z, \omega) \rightarrow s(z) \text { for all } \omega \in N^c \text { and } z \in \mathbb{C}_m^{+}.
$$
Since the convergence above holds for every $m$, we conclude that
$$
s_n(z, \omega) \rightarrow s(z) \text { for all } \omega \in N^c \text { and } z \in \mathbb{C}^{+}.
$$
Thus, for all $z \in \mathbb{C}^{+}$, $s_n^{\bB_n}(z) \rightarrow s(z)$ almost surely. By Theorem B.9 in \cite{BSbook}, we conclude that
$$
F^{\bB_n} \xrightarrow{w} F \text {, a.s. }.
$$
Thus we complete the proof of Theorem \ref{thm:LSD}.

\subsection{Proof of Theorem \ref{cor-extreigen}}\label{se:proof-largest}
Since $\lambda_1^{\mathbf{B}_n}=\sqrt{\frac{N}{p}} \lambda_1^{\mathbf{R}_n}-\sqrt{\frac{p}{N}},$ Theorem \ref{cor-extreigen} can be obtained directly by Lemma 1 and 7 in \cite{gao2017high} when $p/n\to c\in(0,\infty)$. Then we focus on the case where $p / n \rightarrow \infty$. 

\textbf{Proof of Theorem \ref{cor-extreigen} (i) :}\\
By Theorem \ref{thm:LSD}, we have $\liminf _{n\rightarrow \infty} \lambda_{1 }\left(\mathbf{B}_n\right) \geq 1$ a.s..
% $$
% \liminf _{n\rightarrow \infty} \lambda_{1 }\left(\mathbf{B}_n\right) \geq 1 \quad \text { a.s..}
% $$
Thus to prove conclusion (i) in 
Theorem \ref{cor-extreigen}, it suffices to show that
$$
\limsup_{n\rightarrow \infty}  \lambda_{1 }\left(\mathbf{B}_n\right) \leq 1
\quad \text { a.s..}
$$
Firstly, according to Assumption \ref{asp1*}, we truncate the underlying random variables. Here, we choose $\delta_n$
satisfying 
\begin{align}\label{trunc-delta}
\delta_n^{-2(t+1)} \mathbb{E}\left|x_{11}\right|^{2t+2} \mathbbm{1}_{\left\{\left|x_{11}\right| \geqslant \delta_n (n p)^{1/(2t+2)}\right\}} \rightarrow 0, \quad \delta_n  \rightarrow  0, \quad \delta_n (n p)^{1/(2t+2)}  \rightarrow \infty.\end{align}
Similar as arguments in section \ref{se:proof-lsd}, 
let $\hat\bB_n=\hat\bB_n(\hat x_{ij})$,  
$\tilde\bB_n=\tilde\bB_n(\tilde x_{ij})$
be defined similarly to $\bB_n$ with $x_{ij}$ replaced by $\hat{x}_{ij}$, and 
$\tilde{x}_{ij}$ respectively,  where $\hat x_{ij}=x_{ij}I_{\{|x_{ij}|\leq \delta_n (n p)^{1/(2t+2)}\}}$
and $\tilde{x}_{ij}=(\hat x_{ij}-\E\hat x_{ij})/\sigma_n$ with $\sigma_n^2=\E|\hat x_{ij}-\E\hat x_{ij}|^2$. Following the proof of Theorem 1 in \cite{chenbingbing}, we have 
$$
\mathbb{P}\left(\mathbf{B}_n \neq \hat{\mathbf{B}}_n, \text { i.o. }\right)=0\quad \text { a.s.},
$$
from which we obtain $\lambda_{1}\left(\mathbf{B}_n\right)-\lambda_{1 }\left(\hat{\mathbf{B}}_n\right) \rightarrow 0$ a.s. as $n\rightarrow \infty$. And note that $\hat\bB_n=\tilde\bB_n$. We have $\lambda_{1}\left(\mathbf{B}_n\right)-\lambda_{1 }\left(\tilde{\mathbf{B}}_n\right) \rightarrow 0$  a.s.. By the above results, it is sufficient
to show that $\limsup \sup _{n \rightarrow \infty} \lambda_{1}\left(\tilde{\mathbf{B}}_n\right) \leq 1$ a.s.. To this end, note that $\tilde\bB_n$ satisfies truncation condition of Theorem \ref{cor-extreigen} (ii). Therefore, we can obtain Theorem \ref{cor-extreigen} (i) according to conclusion (ii). Next we give the proof of  the conclusion (ii).

\textbf{Proof of Theorem \ref{cor-extreigen} (ii):}\\
To begin with, by (S11) in \cite{yulong}, we have
$$
\lambda_1^{\bB_n}\leq \lambda_1^{\bPhi^2}\lambda_1^{\bB_{n0}}\leq \max_{1\leq i\leq n}\left|
\be_i^{\top}\bB_{n0}\be_i
\right|+\lambda_1^{\bC_n},
$$
where $$\mathbf{B}_{n 0}=\sqrt{\frac{p}{N}}\left[\frac{1}{p} \mathbf{X}^{\top} \mathbf{D}_n \mathbf{X}-\mathbf{I}_n\right],\quad
\bC_n=\bB_{n0}-\operatorname{diag}(\bB_{n0}).
$$
Since $$\mathbf{e}_i^{\top} \mathbf{B}_{n 0} \mathbf{e}_i=
\frac{1}{\sqrt{pN}}\sum_{k=1}^p(\frac{1}{s_{kk}}X_{ki}^2-1)=
\frac{1}{\sqrt{pN}}\sum_{k=1}^p\frac{1}{s_{kk}}(X_{ki}^2-1)+
\frac{1}{\sqrt{pN}}\sum_{k=1}^p(\frac{1}{s_{kk}}-1)
, 
$$
to prove Theorem \ref{cor-extreigen}, it is sufficient to prove, for any $\epsilon>0, \ell>0$,
\begin{align}
&\mathbb{P}\left(\max _{1\leq i \leq n} \frac{1}{\sqrt{p N}} \sum_{k=1}^p \left|\frac{1}{s_{k k}}\left(X_{k i}^2-1\right)\right|>\epsilon\right)=\mathrm{o}\left(n^{-\ell}\right),
\label{n-ell-1}
\\
&\mathbb{P}\left(\frac{1}{\sqrt{p N}} \sum_{k=1}^p \left|\frac{1}{s_{k k}}-1\right|>\epsilon\right)=\mathrm{o}\left(n^{-\ell}\right),
\label{n-ell-2}\\
&\mathbb{P}\left(\lambda_1^{\bC_n}>2+\epsilon\right)=\mathrm{o}\left(n^{-\ell}\right).
\label{n-ell-3}
\end{align}

The proofs of \eqref{n-ell-1}-\eqref{n-ell-3} rely on Lemma
\ref{le-max-eig} below. The proof of Lemma \ref{le-max-eig} is postponed to the supplementary material.
\begin{lemma}\label{le-max-eig} Under the assumptions of Theorem \ref{cor-extreigen} (ii), we have
    \begin{align*}
\mathbb{P}\left(\max_{1\leq k\leq p}\left|\frac{1}{s_{kk}}-1\right|>\epsilon\right)=o(n^{-\ell}).
\end{align*}
\end{lemma}

By Lemma \ref{le-max-eig}, $\max_{1\leq k\leq p} 1/s_{kk}<2$ with high probability, then \eqref{n-ell-1} comes directly from (9) in \cite{chenbingbing}. For \eqref{n-ell-2}, by Burkholder inequality (Lemma 2.13 in \cite{BSbook}), we have
\begin{align*}
&\mathbb{P}\left(\frac{1}{\sqrt{p N}} \sum_{k=1}^p\left|\frac{1}{s_{k k}}-1\right|>\epsilon\right)=
\mathbb{P}\left(\sum_{k=1}^p\left|\frac{1}{s_{k k}}-1\right|>\epsilon\sqrt{p N}\right)\\
&\leq 
K\frac{\E\left|\sum_{k=1}^p(s_{kk}-1)\right|^\ell}{(\epsilon\sqrt{p N})^\ell}+o(n^{-\ell})\\
&\leq 
K\frac{
\left(\sum_{k=1}^p\E\left|s_{kk}-1\right|^2\right)^{\ell/2}
+\sum_{k=1}^p\E\left|s_{kk}-1\right|^\ell
}{(\epsilon\sqrt{p N})^\ell}+o(n^{-\ell})\\
&\leq
K\frac{
(p/n)^{\ell/2}+pn^{-\ell/2}+pn^{-\ell+1}v_{2\ell}
}{(\epsilon\sqrt{p N})^\ell}+o(n^{-\ell})=o(n^{-\ell}).
\end{align*}
And for \eqref{n-ell-3}, by using Lemma \ref{le-max-eig} again, we have
for any $\epsilon,\epsilon^{\prime}>0$,
\begin{align*}
\mathbb{P}\left(\lambda_1^{\mathbf{C}_n}>2+\epsilon\right)&=
\mathbb{P}\left(\lambda_1^{\mathbf{C}_n}>2+\epsilon,\max _{1 \leq k \leq p}\left|\frac{1}{s_{k k}}-1\right|<\epsilon^{\prime}\right)+
\mathbb{P}\left(\lambda_1^{\mathbf{C}_n}>2+\epsilon,\max _{1 \leq k \leq p}\left|\frac{1}{s_{k k}}-1\right|>\epsilon^{\prime}\right)\\
&=\mathbb{P}\left(\lambda_1^{\mathbf{C}_n}>2+\epsilon,\max _{1 \leq k \leq p}\left|\frac{1}{s_{k k}}-1\right|<\epsilon^{\prime}\right)+o\left(n^{-\ell}\right)\\
&=o\left(n^{-\ell}\right),
\end{align*}
where the last equality holds by (8) in \cite{chenbingbing} and (S12) in \cite{yulong}.  Together with \eqref{n-ell-1} and \eqref{n-ell-2}, we obtain $\Prob\left(\lambda_1\left(\mathbf{B}_n\right) \geq 2+\epsilon\right)=\mathrm{o}\left(n^{-\ell}\right)$.
Therefore we complete the proof. 
%of Theorem \ref{cor-extreigen}.

\subsection{Proof of Theorem \ref{ult-cltmain}}

Now we present the strategy for the proof of Theorem \ref{ult-cltmain}. By the Cauchy integral formula, we have
$
\int f(x)\dif G(x)=-\frac{1}{2\pi i}\oint_\mathcal{C}f(z)m_G(z)\dif z
$
valid for any c.d.f $G$ and any analytic function $f$ on an open set containing the support of $G$, where $\oint_\mathcal{C}$ is the contour integration in the anti-clockwise direction. In our case, $G(x)=n(F^{\bB_n}(x)-F^{c_N}(x))$.  Therefore, the problem of finding the limiting distribution reduces to the study of $M_n^{\bB_n}(z)$:
\begin{align*}
    M_n^{\bB_n}(z)=n\left(s_n^{\bB_n}(z)-s_{c_N}(z)\right)-\Theta_n (s_{c_N}(z)).
\end{align*}
%Since  
%\begin{align*}
%s_n^{\bB_n}(z)=s_n(z)-\frac{1}{n}\frac{\sqrt{c_N}}{z(\sqrt{c_N}+z)},
%\end{align*}
%the CLT for LSS of $\bB_n$ can be obtained  by the main results in \cite{gao2017high} when $p/n\to c\in(0,\infty)$. We thus focus on the case where $p/n\to \infty$.

By using \eqref{sb-sn}, under the ultrahigh dimensional case,
$$
\Theta_n (s_{c_N}(z))=\frac{s^3(z)+s(z)-s^{\prime}(z) s(z)}{s^2(z)-1}-\frac{1}{z}+o(1),
$$
then we have
\begin{align}\label{MBn}
M_n^{\bB_n}(z)=M_n(z)-
\frac{s^3(z)+s(z)-s^{\prime}(z) s(z)}{s^2(z)-1}+o(1),
\end{align}
where
\begin{align*}
    M_n(z) = n\left(s_n(z)-s_{c_N}(z)\right).
\end{align*}
Firstly, according to Assumption \ref{asp1*}, we truncate the underlying random variables. Here, we choose $\delta_n$
defined in \eqref{trunc-delta}.
By the  arguments in section \ref{se:proof-largest}, 
let $\hat\bB_n=\hat\bB_n(\hat x_{ij})$,  
$\tilde\bB_n=\tilde\bB_n(\tilde x_{ij})$
be defined similarly to $\bB_n$ with $x_{ij}$ replaced by $\hat{x}_{ij}$, and 
$\tilde{x}_{ij}$ respectively,  where $\hat x_{ij}=x_{ij}I_{\{|x_{ij}|\leq \delta_n (n p)^{1/(2t+2)}\}}$
and $\tilde{x}_{ij}=(\hat x_{ij}-\E\hat x_{ij})/\sigma_n$ with $\sigma_n^2=\E|\hat x_{ij}-\E\hat x_{ij}|^2$.
We then conclude that
$$
P\left(\mathbf{B}_n \neq \hat{\mathbf{B}}_n\right) \leq n p \cdot P\left(\left|x_{11}\right| \geq \delta_n (n p)^{1/(2t+2)}\right) \leq K \delta_n^{-2(t+1)} \mathbb{E}\left|x_{11}\right|^{2 t+2} \mathbbm{1}_{\left\{\left|x_{11}\right| \geqslant \delta_n(n p)^{1 /(2 t+2)}\right\}}=o(1).
$$
Let $\hat{G}_n(f)$ and $\tilde{G}_n(f)$ be $G_n(f)$ with $\mathbf{B}_n$ replaced by $\hat{\mathbf{B}}_n$ and $\tilde{\mathbf{B}}_n$ respectively. Then for each $j=1,2, \ldots, k,$ since $\hat\bB_n=\tilde\bB_n$, we have
$$
G_n(f_j)=\hat{G}_n(f_j)+o_p(1)=\tilde{G}_n(f_j)+o_p(1).
$$
Thus, we only need to find the limit distribution of $\left\{\widetilde{G}_n\left(f_j\right), j=1, \ldots, k\right\}$. Hence, in what follows, we assume that the underlying variables are truncated at $\delta_n (n p)^{\frac{1}{2t+2}}$, centralized, and renormalized. For convenience, we shall suppress the superscript on the variables, and assume that, for any $1 \leqslant i \leqslant p$ and $1 \leqslant j \leqslant n$,
\begin{align}\label{con:trunc}
& \left|x_{i j}\right| \leqslant \delta_n (n p)^{\frac{1}{2t+2}}, \quad \mathbb{E} x_{i j}=0, \quad \mathbb{E} |x_{i j}|^2=1, \quad
\E |x_{i j}|^4=\kappa+o(1),\quad
\E |x_{i j}|^{2t+2}<\infty.
\end{align}

For any $\varepsilon>0$, define the event $F_n(\varepsilon)=\left\{\max _{j \leq n}\left|\lambda_j(\mathbf{B}_n)\right| \geq 2+\varepsilon\right\}$ where $\mathbf{B}_n$ is defined by the truncated and normalized variables satisfying condition \eqref{con:trunc}. By Theorem \ref{cor-extreigen}, for any $\ell>0$
\begin{align}\label{n-ell}
\Prob\left(F_n(\varepsilon)\right)=\mathrm{o}\left(n^{-\ell}\right) .
\end{align}
Here we would point out that the result regarding the minimum eigenvalue of $\mathbf{B}_n$ can be obtained similarly by investigating the maximum eigenvalue of $-\mathbf{B}_n$.

Note that the support of $F^{\bB_{n}}$ is random. Fortunately, we have shown that the extreme eigenvalues of $\bB_{n}$ are highly concentrated around two edges of the support of the limiting semicircle law $F(x)$ in \eqref{n-ell}. Then the contour $\mathcal{C}$ can be appropriately chosen. Moreover,  as in \cite{bai2004clt}, by \eqref{n-ell}, we can replace the process $\{M_n(z),z\in\mathcal{C}\}$ by a slightly modified process $\{\widehat{M}_n(z), z\in\mathcal{C}\}$. Below we present the definitions of the contour $\mathcal{C}$ and the modified process $\widehat{M}_n(z)$. Let $x_r$ be any number greater than $2+\frac{1}{\sqrt{c_N}}$. Let $x_l$ be any number less than $-2+\frac{1}{\sqrt{c_N}}$. Now let $\mathcal{C}_u=\{x+iv_0:x\in[x_l,x_r]\}$. Then we define $\mathcal{C}^+ := \{x_l+iv:v\in[0,v_0]\}\cup\mathcal{C}_u\cup\{x_r+iv: v\in[0,v_0]\}$, and $\mathcal{C}=\mathcal{C}^+\cup \overline{\mathcal{C}^+}$. Now we define the subsets $\mathcal{C}_n$ of $\mathcal{C}$ on which $M_n(\cdot)$ equals to $\widehat{M}_n(\cdot)$. Choose sequence $\{\varepsilon_n\}$ decreasing to zero satisfying for some $\alpha\in(0,1)$, $\varepsilon_n\geq n^{-\alpha}$. Let
\begin{equation*}
\mathcal{C}_l=\{x_l+iv:v\in[0,v_0]\},
\end{equation*}
% \begin{equation*}
% \mathcal{C}_l=\left\{
% \begin{array}{ccc}
% \{x_l+iv:v\in[n^{-1}\varepsilon_n,v_0]\}~~~~ \mathrm{if}~x_l>0,\\
% \{x_l+iv:v\in[0,v_0]\} ~~~~~~~~~~\mathrm{if}~x_l<0,
% \end{array}
% \right.
% \end{equation*}
and $\mathcal{C}_r=\{x_r+iv: v\in[n^{-1}\varepsilon,v_0]\}$. Then $\mathcal{C}_n=\mathcal{C}_l\cup\mathcal{C}_u\cup\mathcal{C}_r$. For $z=x+iv$, we define
% \begin{equation*}
% \widehat{M}_p(z)=\left\{
% \begin{array}{ccc}
% &M_p(z)\mathrm{for}~ z\in \mathcal{C}_p,\\
% &M_p(x_r+in^{-1}\varepsilon_n)\ \ \ \ \mathrm{for}~ x=x_r,v\in[0,n^{-1}\varepsilon_n], and \ if x_l>0,\\
% &M_p(x_l+in^{-1}\varepsilon_n)\mathrm{for}~
% x=x_l,v\in[0,n^{-1}\varepsilon_n].
% \end{array}\right.
% \end{equation*}
\begin{equation*}
\widehat{M}_n(z)= \begin{cases}M_n(z), & \text { for } z \in \mathcal{C}_n, \\ 
M_n(x_r+i n^{-1} \varepsilon_n), & \text { for } x=x_r, v \in[0, n^{-1} \varepsilon_n],\\ 
M_n(x_l+i n^{-1} \varepsilon_n), & \text { for } x=x_l,v\in[0,n^{-1}\varepsilon_n].\end{cases}
\end{equation*}
With the help of \eqref{n-ell}, one may thus find 
$$
\oint_{\mathcal{C}} f_{j}(z) M_n(z) d z=\oint_{\mathcal{C}} f_{j}(z) \widehat{M}_n(z) d z+o_p(1),
$$
for every $j\in\{1,\dots,K\}$. Hence according to \eqref{MBn}, the proof of Theorem \ref{ult-cltmain} can be completed by verifying the convergence of $\widehat{M}_n(z)$ on ${\mathcal{C}}$ as stated in the following lemma.

% Let
% \begin{align*}
% M_n^{\bB_n}(z)=n\left(s_n^{\bB_n}(z)-s_{c_N}(z)\right)-\Theta_n (s_{c_N}(z)).
% \end{align*}

\begin{lemma}\label{cormpzclt}
     In addition to Assumptions \ref{asp1*}, \ref{asp2},\ref{asp3}, 
suppose condition \eqref{con:trunc} holds.
We have 
$$
\widehat{M}_n(z) \stackrel{d}{=} M(z)+o_p(1), \quad z \in \mathcal{C},
$$
where the random process $M(z)$ is a two-dimensional Gaussian process. The mean function is
$$
\E M(z)=\frac{s^3(z)+s(z)-s^{\prime}(z) s(z)}{s^2(z)-1},
$$
and the covariance function is
\begin{align}\label{covreal}
    Cov(M(z_1),M(z_2))= 2 \left[\frac{s^{\prime}(z_1)s^{\prime}(z_2)}{\{s(z_1)-s(z_2)\}^2}-\frac{1}{(z_1-z_2)^2}\right] - 2s^{\prime}(z_1)s^{\prime}(z_2).
\end{align}

%after truncation of the
%data, the empirical process $\{M_n^{\bB_n}(z), z \in \mathcal{C}\}$ converges weakly to a centred Gaussian process
%$\{M(z), z \in \mathcal{C}\}$ with the covariance function
%\begin{align}\label{covreal}
   % Cov(M(z_1),M(z_2))= 2 \left[\frac{s^{\prime}(z_1)s^{\prime}(z_2)}{\{s(z_1)-s(z_2)\}^2}-\frac{1}{(z_1-z_2)^2}\right] - 2s^{\prime}(z_1)s^{\prime}(z_2).
%\end{align}
\iffalse
\begin{align}\label{covreal}
    Cov(M(z_1),M(z_2))= 2 \left[\frac{s^{'}(z_1)s^{'}(z_2)}{\{s(z_1)-s(z_2)\}^2}-\frac{1}{(z_1-z_2)^2}\right] - (|\Psi|^2+1)s^{'}(z_1)s^{'}(z_2).
\end{align}
under the real random-variable
case, where $\Psi=\frac{\Expe \{y_{i1}-\Expe y_{i1}\}^2}{\Expe |y_{i1}-\Expe y_{i1}|^2}=1$,
$s(z)$ is the Stieltjes transform of the semicircle law. When ${y_{ij}}$ are complex variables, assuming that $\Psi$ are the same for
$i=1,2, \ldots, p$, the covariance function is
\begin{align}\label{covcomp}
    Cov(M(z_1),M(z_2))&= \eqref{covreal}-[\frac{s^{'}(z_1)s^{'}(z_2)}{\{s(z_1)-s(z_2)\}^2}-\frac{1}{(z_1-z_2)^2}]\nonumber\\
    &+|\Psi|^2\frac{s^2(z_1)s^2(z_2)}{(-1+s^2(z_1))(-1+s^2(z_2))(-1+|\Psi|^2s(z_1)s(z_2))^2}.
\end{align}
\fi
\end{lemma}

To prove Lemma \ref{cormpzclt}, we decompose $\widehat{M}_n(z)$ into a random part $M^{(1)}_n(z)$ and a deterministic part $M^{(2)}_n(z)$ for $z\in\mathcal{C}_n$, that is, $M_n(z)=M_n^{(1)}(z)+M_n^{(2)}(z)$, where  
\[
M_n^{(1)}(z)=n\big\{s_{n}(z)-\mathbb{E}s_{n}(z)\big\}
\quad\text{and}\quad
M_n^{(2)}(z)=n\big\{\mathbb{E}s_{n}(z)-s_{c_N}(z)\big\}.
\]

The random part contributes to the covariance function and the deterministic part contributes
to the mean function. By Theorem 8.1 in \citet{billingsley1968convergence}, the proof of Lemma \ref{cormpzclt} is
then complete if we can verify the following three steps:
\begin{itemize}
	%\item[Step 1] Truncation, centralization, and rescaling.
	\item[\bf Step 1.] Finite-dimensional convergence of $M_n^{(1)}(z)$ in distribution on $\mathcal{C}_n$  to a centered multivariate Gaussian random vector with covariance function given by \eqref{covreal}.
    \begin{lemma}\label{step1}
    Under assumptions of Theorem \ref{ult-cltmain} and condition \eqref{con:trunc}, as $n\rightarrow\infty$, for any set of $r$ points $\{z_1, z_2, . . ., z_r\} \subseteq \mathcal{C}$, the random vector
$\bigl(M_n^{(1)}(z_1), \ldots ,M_n^{(1)}(z_r)\bigr)$ converges weakly to a $r$-dimensional centered Gaussian distribution with covariance function in \eqref{covreal}.
\end{lemma}
 %or \eqref{covcomp}.
	\item[\bf Step 2.] Tightness of the $M_n^{(1)}(z)$ for $z\in\mathcal{C}_n$. The tightness can be established by Theorem 12.3 of \cite{billingsley1968convergence}. It's sufficient to 
    verify the moment condition given in the following lemma.
\begin{lemma}\label{step2}
    Under assumptions of Lemma \ref{step1},
$
\sup _{n ; z_1, z_2 \in \mathcal{C}_n} \frac{\mathbb{E}\left|M_n^{(1)}\left(z_1\right)-M_n^{(1)}\left(z_2\right)\right|^2}{\left|z_1-z_2\right|^2}<\infty.
$
\end{lemma}   
	\item[\bf Step 3.] Convergence of the non-random part $M_n^{(2)}(z)$.
    \begin{lemma}\label{step3}
      Under assumptions of Lemma \ref{step1},
 $M_n^{(2)}(z)=\frac{s^3(z)+s(z)-s^{\prime}(z) s(z)}{s^2(z)-1}+o(1)$ for $z\in\mathcal{C}_n$.  
 %$M_n^{(2)}(z)=\Theta_n (s_{c_N}(z))+\frac{\sqrt{c_N}}{z(\sqrt{c_N}+z)}+o(1)$ for $z\in\mathcal{C}_n$. 
    \end{lemma}
\end{itemize}
Thus we complete the proof of Theorem \ref{ult-cltmain}. The proof of Lemma~\ref{step1} is presented below while the proofs of Lemma~\ref{step2}-\ref{step3} are delegated to the supplement file due to page limit.

\subsection{Proof of Lemma~\ref{step1}}

To prove Lemma \ref{step1}, we first introduce the following supporting lemmas.
\begin{lemma}\label{le-supp}
    Under assumptions of Lemma \ref{step1}, we have
    \begin{align*}
\mathcal{Y}_1\left(z_1, z_2\right)&\triangleq -\frac{\partial^2}{\partial z_1\partial z_2}\Big(\sum^{p}_{j=1}\big[\Expe_{j-1}\tilde{\beta}_j(z_1)\varepsilon_j(z_1)][\Expe_{j-1}\tilde{\beta}_j(z_2)\varepsilon_j(z_2)\big]\Big)=o_p(1),\\
\mathcal{Y}_2\left(z_1,z_2\right)&\triangleq\frac{\partial^2}{\partial z_1\partial z_2}\Big(\sum^{p}_{j=1}\Expe_{j-1}\big[\Expe_j\big(\tilde{\beta}_j(z_1)\varepsilon_j(z_1)\big)\Expe_j\big(\tilde{\beta}_j(z_2)\varepsilon_j(z_2)\big)\big]\Big)\\
&=2\frac{\partial^2}{\partial z_1 \partial z_2} \mathcal{J}-2s^{\prime}\left(z_1\right) s^{\prime}\left(z_2\right)+o_p(1),
    \end{align*}
    where$$
\mathcal{J} =\frac{1}{n^2} b_n\left(z_1\right) b_n\left(z_2\right)\left[\mathbb{E} \sum_{j=1}^p \operatorname{tr}\left[\mathbb{E}_j\left(\mathbf{Q}_j^{-1}\left(z_1\right)\right) \mathbb{E}_j\left(\mathbf{Q}_j^{-1}\left(z_2\right)\right)\right]\right].
    $$
\end{lemma}

\begin{lemma}\label{le-J}
     Under assumptions of Lemma \ref{step1}, we have
     $$\frac{\partial^2}{\partial z_2 \partial z_1} \mathcal{J}=\frac{s^2\left(z_1\right) s^2\left(z_2\right)}{\left[s^2\left(z_1\right)-1\right]\left[s^2\left(z_2\right)-1\right]\left[s\left(z_1\right) s\left(z_2\right)-1\right]^2}+o_p(1).$$
\end{lemma}

The proof of Lemma~\ref{le-J} is presented in next section while the proof of Lemma~\ref{le-supp} is delegated to the supplement file due to page limit.

We now proceed to the proof of Lemma~\ref{step1}. By the fact that a random vector is multivariate normally distributed if and only if every linear combination of its components is normally distributed, we need only show that, for any positive integer $r$ and any complex sequence
${a_j}$, the sum 
\[
    \sum_{j=1}^ra_jM^{(1)}_n(z_j)
\]
converges weakly to a Gaussian random variable. To this end, we first decompose the random part $M^{(1)}_n(z)$ as a sum of martingale difference, which is given in \eqref{mcd1}. Then, we apply the martingale CLT (Proposition \ref{MCLT}) to obtain the asymptotic distribution of $M^{(1)}_n(z)$. %Details  are provided in the following.
\begin{Proposition}\label{MCLT}
    (Theorem 35.12 of \cite{billingsley1968convergence}). Suppose for each $n, Y_{n, 1}, Y_{n, 2}$, $\ldots, Y_{n, r_n}$ is a real martingale difference sequence with respect to the increasing $\sigma$-field $\left\{\mathcal{F}_{n, j}\right\}$ having second moments. If as $n \rightarrow \infty$, 
(i) $\sum_{j=1}^{r_n} E\left(Y_{n, j}^2 \mid \mathcal{F}_{n, j-1}\right) \xrightarrow{i . p .} \sigma^2,$ and (ii) $\sum_{j=1}^{r_n} E\left(Y_{n, j}^2 I_{\left(\left|Y_{n, j}\right| \geq \varepsilon\right)}\right) \rightarrow 0,$
where $\sigma^2$ is a positive constant and $\varepsilon$ is an arbitrary positive number, then
$
\sum_{j=1}^{r_n} Y_{n, j} \xrightarrow{D} N\left(0, \sigma^2\right) .
$
\end{Proposition}

To begin with, similar as \eqref{estimate-lsd}, we give some useful estimate below.
For $q>2$, we have
\begin{align}\label{estimate}
&\E\left|\varepsilon_k(z)\right|^q\leq K\left(n^{-q/2}+n^{-q/2} p^{q/2-1} \delta_n^{2q-4}\right),\
\E\left|\delta_k(z)\right|^q\leq K\left(n^{-q/2}+n^{-q/2} p^{q/2-1} \delta_n^{2q-4}\right),
\nonumber\\
&\E\left|\tilde\beta_k(z)-b_n(z)\right|^q=O(n^{q/2}p^{-q}),\
\left|b_n(z)-b_1(z)\right|=O(n^{1/2}p^{-2/3}),\
\E\left|b_n(z)-\E\beta_k(z)\right|=O(np^{-2}),\nonumber\\
&\E\left|\gamma_k(z)-\varepsilon_k(z)\right|^q=O(n^{-q/2}),\ \E\left|\frac{1}{n} \operatorname{tr}\left(\bQ^{-1}(z)\bM\right)-\E \frac{1}{n} \operatorname{tr}\left(\bQ^{-1}(z) \bM\right)\right|^q=O(-n^{q/2}).
\end{align}
 Write $M_n^{(1)}(z)$ as a sum of martingale difference sequences (MDS), and then utilize the CLT of MDS  to derive the asymptotic distribution of $M_n^{(1)}(z)$, which can be written as
\begin{align}\label{m1}
M_n^{(1)}(z)&=n[s_{n}(z)-\E s_{n}(z)]
= \sum^{p}_{j=1}(\Expe_j-\Expe_{j-1})\tr\left[\bQ^{-1}(z)-\bQ^{-1}_j(z)\right]\nonumber\\
&= -\sum^{p}_{j=1}(\Expe_j-\Expe_{j-1})\beta_j(z)\br_j^{\top}\bQ_j^{-2}(z)\br_j.
\end{align}
By using \eqref{expand-beta} and the fact that $\left(\mathbb{E}_j-\mathbb{E}_{j-1}\right) \tilde{\beta}_j(z) \frac{1}{n} \operatorname{tr} \mathbf{Q}_j^{-2}(z)=0$, we have 
\begin{align*}
   & \left(\mathbb{E}_j-\mathbb{E}_{j-1}\right) \beta_j(z) \mathbf{r}_j^{\top} \mathbf{Q}_j^{-2}(z) \mathbf{r}_j
 =\mathbb{E}_j\left[\tilde{\beta}_j(z) \delta_j(z)-\tilde{\beta}_j^2(z) \varepsilon_j(z) \frac{1}{n} \operatorname{tr} \mathbf{Q}_j^{-2}(z)\right]\\&+\E_{j-1}[Y_j(z)] 
 -\left(\mathbb{E}_j-\mathbb{E}_{j-1}\right)\left[\tilde{\beta}_j^2(z)\left(\varepsilon_j(z) \delta_j(z)-\beta_j(z) \mathbf{r}_j^{\top}\mathbf{Q}_j^{-2}(z) \mathbf{r}_j \varepsilon_j^2(z)\right)\right],
\end{align*}
where $Y_j(z)=-\mathbb{E}_j\left[\tilde{\beta}_j(z) \delta_j(z)-\tilde{\beta}_j^2(z) \varepsilon_j(z) \frac{1}{n} \operatorname{tr} \mathbf{Q}_j^{-2}(z)\right]$.
With the help of \eqref{estimate}, we have
$$
\begin{aligned}
\mathbb{E}\left|\sum_{j=1}^p\left(\mathbb{E}_j-\mathbb{E}_{j-1}\right) \widetilde{\beta}_j^2(z) \varepsilon_j(z) \delta_j(z)\right|^2 & =\sum_{j=1}^p \mathbb{E}\left|\left(\mathbb{E}_j-\mathbb{E}_{j-1}\right) \widetilde{\beta}_j^2(z) \varepsilon_j(z) \delta_j(z)\right|^2 \\
& \leq K \sum_{j=1}^p \mathbb{E}\left|\widetilde{\beta}_j^2(z) \varepsilon_j(z) \delta_j(z)\right|^2
\leq K(p^{-1}+\delta_n^4)
=o(1),
\end{aligned}
$$
and similarly
$$
\mathbb{E}\left|\sum_{j=1}^p\left(\mathbb{E}_j-\mathbb{E}_{j-1}\right) \widetilde{\beta}_j^2(z) \beta_j(z) \mathbf{r}_j^{\top} \mathbf{D}_j^{-2}(z) \mathbf{r}_j \varepsilon_j^2(z)\right|^2=o(1).
$$
By \eqref{m1}, we obtain 
\begin{align}\label{mcd1}
    M_n^{(1)}(z)=\sum_{j=1}^p[\E_{j}-\E_{j-1}]Y_j(z)+o_p(1).
\end{align}
Then we need to consider the limit of the following term:
$$
\sum_{i=1}^r \alpha_i \sum_{j=1}^p \left[\mathbb{E}_j-\mathbb{E}_{j-1}\right] Y_j(z)=\sum_{j=1}^p \sum_{i=1}^r \alpha_i \left[\mathbb{E}_j-\mathbb{E}_{j-1}\right] Y_j(z).
$$
Using \eqref{estimate} we obtain
$$
\mathbb{E}\left|Y_j(z)\right|^4 \leq K \left(
c_n^{-2}\mathbb{E}\left|\delta_j(z)\right|^4
+c_n^{-4}\mathbb{E}\left|\varepsilon_j(z)\right|^4\right)
\leq K\left(p^{-2}+p^{-1}\delta_n^4\right),
$$
from which we can have
\begin{align*}
&\sum_{j=1}^p \mathbb{E}\left(\left|\sum_{i=1}^r \alpha_i \left[\mathbb{E}_j-\mathbb{E}_{j-1}\right] Y_j(z_i)\right|^2 I_{\left(\left|\sum_{i=1}^r \alpha_i \left[\mathbb{E}_j-\mathbb{E}_{j-1}\right] Y_j(z_i)\right| \geq \varepsilon\right)}\right) \\
&\leq \frac{1}{\varepsilon^2} \sum_{j=1}^p \mathbb{E}\left|\sum_{i=1}^r \alpha_i \left[\mathbb{E}_j-\mathbb{E}_{j-1}\right] Y_j(z_i)\right|^4\to 0.
\end{align*}
Thus the condition (ii) of Proposition \ref{MCLT} is satisfied.

Next, it suffices to prove that
\begin{align}\label{condi2eq}
\sum_{j=1}^p \mathbb{E}_{j-1}\left[Y_j\left(z_1\right)-E_{j-1}Y_j\left(z_1\right)\right]
\left[Y_j\left(z_2\right)-E_{j-1}Y_j\left(z_2\right)\right]
\end{align}
converges in probability to \eqref{covreal}.
Note that
\begin{align*}
    \eqref{condi2eq}=\sum^{p}_{j=1}\Expe_{j-1}[Y_j(z_1)Y_j(z_2)]-\sum^{p}_{j=1}[\Expe_{j-1}Y_j(z_1)][\Expe_{j-1}Y_j(z_2)]
    =\mathcal{Y}_1\left(z_1, z_2\right)+\mathcal{Y}_2\left(z_1, z_2\right).
\end{align*}
By Lemmas \ref{le-supp}-\ref{le-J}, we obtain the limit of \eqref{condi2eq} is \eqref{covreal}. Thus we complete the proof of Lemma \ref{step1}.

\iffalse
where 
\begin{align*}
\mathcal{Y}_1\left(z_1, z_2\right)&= -\frac{\partial^2}{\partial z_1\partial z_2}\Big(\sum^{p}_{j=1}\big[\Expe_{j-1}\tilde{\beta}_j(z_1)\varepsilon_j(z_1)][\Expe_{j-1}\tilde{\beta}_j(z_2)\varepsilon_j(z_2)\big]\Big),\\
\mathcal{Y}_2\left(z_1, z_2\right)&=\frac{\partial^2}{\partial z_1\partial z_2}\Big(\sum^{p}_{j=1}\Expe_{j-1}\big[\Expe_j\big(\tilde{\beta}_j(z_1)\varepsilon_j(z_1)\big)\Expe_j\big(\tilde{\beta}_j(z_2)\varepsilon_j(z_2)\big)\big]\Big).
\end{align*}
Thus, it's enough to consider the limits of $\mathcal{Y}_i\left(z_1, z_2\right), i=1,2$, as given in the following lemma.

By Lemmas \ref{le-supp}-\ref{le-J}, we obtain the limit of \eqref{condi2eq}. Thus we complete the proof of Lemma \ref{step1}.

The proof of Lemma \ref{le-supp} is postponed to the supplementary material. Considering $\mathcal{J}$ next. In contrast to high dimensional case where $p/n\to c\in(0,\infty)$ \citep{gao2017high}, we consider the ultra-high dimensional setting with $p/n\to\infty$. In this regime, we carefully analyze the effects brought by $c_n$ and $c_N$, and derive a completely new determinant equivalent form for $\mathbf{Q}_j^{-1}(z)$, the resolvent of the renormalized correlation matrix with the $j$th component information removed.
\fi
\subsection{Proof of Lemma \ref{le-J}}

The proof of Lemma \ref{le-J} differs substantially from the classical case. Unlike the high dimensional case where $p/n \to c \in (0, \infty)$ \citep{gao2017high}, our analysis is conducted in the ultrahigh dimensional regime with $p/n \to \infty$. In this setting, we carefully examine the influence of $c_n$ and $c_N$, and derive a novel determinant equivalent form
$\left(\frac{p-1}{N} b_1\left(z\right)-\sqrt{c_N}-z\right)^{-1} \mathbf{I}_n$
for $\mathbf{Q}_j^{-1}(z)$, the resolvent of the renormalized correlation matrix with the $j$th component information removed.  

Specifically, by using the identity
$
\mathbf{r}_i^{\top} \mathbf{Q}_j^{-1}(z)=\sqrt{c_N}\beta_{ij}(z)\mathbf{r}_i^{\top}\mathbf{Q}_{ij}^{-1}(z)
$, we get
$$
\mathbf{Q}_j^{-1}(z)=-\bH_n(z)+b_1(z_1)\bA(z_1)+\bB(z_1)+\bC(z_1)+\bF(z_1),
$$
where $\bH_n(z_1)=\left(\sqrt{c_{N}}+z_1-\frac{p-1}{N} b_1\left(z_1\right)\right)^{-1}\bI_n$ and 
\begin{align*}
& \bA\left(z_1\right)=\sum_{i \neq j}^p \bH_n\left(z_1\right)\left(\br_i \br_i^{\top}-\frac{1}{n-1} \bPhi\right) \bQ_{i j}^{-1}\left(z_1\right), \\
& \bB\left(z_1\right)=\sum_{i \neq j}^p\left(\beta_{i j}\left(z_1\right)-b_1\left(z_1\right)\right) \bH_n\left(z_1\right) \br_i \br_i^{\top} \bQ_{i j}^{-1}\left(z_1\right), \\
& \bC\left(z_1\right)=-\frac{p-1}{N} b_1\left(z_1\right) \bH_n\left(z_1\right) \bPhi\left(\bQ_j^{-1}\left(z_1\right)-\bQ_{i j}^{-1}\left(z_1\right)\right), \\
& \bF\left(z_1\right)=-\frac{p-1}{nN} b_1\left(z_1\right) \bH_n\left(z_1\right) \mathbf{1}_n \mathbf{1}_n^{\top} \bQ_{j}^{-1}\left(z_1\right).
\end{align*}
We next employ $-\bH_n(z)$ as a suitable approximation to the resolvent matrix $\bQ_{j}^{-1}(z)$, extract the dominant terms contributing to the limiting behavior of $\mathcal{J}$, and demonstrate that the error terms are negligible.
Note that $\|\bH_n(z_1)\|\leq K$ and by Lemma 6 in \cite{gao2017high}, we have $\E\mathbf{r}_i \mathbf{r}_i^{\top}=\frac{1}{n-1} \boldsymbol{\Phi}$. 
For any nonrandom $\mathbf{M}$ with $\|\mathbf{M}\|\leq K$, 
by using \eqref{estimate}, we can obtain$$
n^{-1} \E|\tr\bB(z_1)\bM|=O(n^{-1/2}),~
n^{-1} \E|\tr\bC(z_1)\bM|=O(n^{-1}),
$$
which implies 
$$
n^{-1}\E  \left|\tr\mathbb{E}_j\left(\mathbf{B}\left(z_1\right)\right) \mathbf{Q}_j^{-1}\left(z_2\right)\right|=o(1),~n^{-1}\E  \left|\tr\mathbb{E}_j\left(\mathbf{C}\left(z_1\right)\right) \mathbf{Q}_j^{-1}\left(z_2\right)\right|=o(1).
$$
And since $\mathbf{1}_n^{\top} \mathbf{Q}_j^{-1}\left(z_1\right)=-\frac{1}{\sqrt{c_N}+z}\mathbf{1}_n^{\top}$, we have
$
\left|\tr\E_j(\mathbf{F}(z_1))\bQ_j^{-1}(z_2)\right|\leq K/\sqrt{c_n}.
$
In the end, consider the term $
b_1\left(z_1\right)\tr\E_j(\mathbf{A}(z_1))\bQ_j^{-1}(z_2).
$ By using $\mathbf{Q}^{-1}(z)-\mathbf{Q}^{-1}_k(z)=-\mathbf{Q}^{-1}_k(z) \mathbf{r}_k \mathbf{r}_k^{\top} \mathbf{Q}^{-1}_k(z) \beta_k(z)$, we can write  $
\tr\E_j(\mathbf{A}(z_1))\bQ^{-1}_j(z_2)=A_1\left(z_1, z_2\right)+A_2\left(z_1, z_2\right)+A_3\left(z_1, z_2\right),
$ where 
\begin{align*}
A_1\left(z_1, z_2\right) 
& =-\sum_{i<j}^p \beta_{i j}\left(z_2\right) \mathbf{r}_i^{\top} \mathbb{E}_j\left(\mathbf{Q}_{i j}^{-1}\left(z_1\right)\right) \mathbf{Q}_{i j}^{-1}\left(z_2\right) \mathbf{r}_i \mathbf{r}_i^{\top} \mathbf{Q}_{i j}^{-1}\left(z_2\right) \mathbf{H}_n\left(z_1\right) \mathbf{r}_i, \\
A_2\left(z_1, z_2\right) & =-\operatorname{tr} \sum_{i<j}^p \mathbf{H}_n\left(z_1\right) N^{-1} \bPhi\mathbb{E}_j\left(\mathbf{Q}_{i j}^{-1}\left(z_1\right)\right)\left(\mathbf{Q}_j^{-1}\left(z_2\right)-\mathbf{Q}_{i j}^{-1}\left(z_2\right)\right), \\
A_3\left(z_1, z_2\right) & =\operatorname{tr} \sum_{i<j}^p \mathbf{H}_n\left(z_1\right)\left(\mathbf{r}_i \mathbf{r}_i^{\top}-N^{-1}\bPhi \right) \mathbb{E}_j\left(\mathbf{Q}_{i j}^{-1}\left(z_1\right)\right) \mathbf{Q}_{i j}^{-1}\left(z_2\right).
\end{align*}
With \eqref{estimate}, we obtain $\left|b_1(z_1)A_2(z_1,z_2)\right|\leq K$. Our next aim is to show 
\begin{align}\label{A3-o1}
    n^{-1}b_1(z_1)\E_j A_3(z_1,z_2)=o_p(1).
\end{align}
Write \begin{align*}\E\left|b_1(z_1)\E_j A_3(z_1,z_2)\right|^2&=\left|b_1(z_1)\right|^2\sum_{i_1,i_2<j}\E\tr\mathbf{H}_n\left(z_1\right)\left(\mathbf{r}_{i_1} \mathbf{r}_{i_1}^{\top}-N^{-1} \boldsymbol{\Phi}\right) \mathbb{E}_j\left(\mathbf{Q}_{i_1 j}^{-1}\left(z_1\right)\right) \E_j\left(\mathbf{Q}_{i_1 j}^{-1}\left(z_2\right)\right)\\
&\qquad\qquad\qquad\quad\times 
\tr\mathbf{H}_n\left(z_1\right)\left(\mathbf{r}_{i_2} \mathbf{r}_{i_2}^{\top}-N^{-1} \boldsymbol{\Phi}\right) \mathbb{E}_j\left(\mathbf{Q}_{i_2j}^{-1}\left(z_1\right)\right) \E_j\left(\mathbf{Q}_{i_2 j}^{-1}\left(z_2\right)\right)\\
&=\left|b_1(z_1)\right|^2\sum_{i_1,i_2<j}\E\tr\mathbf{H}_n\left(z_1\right)\left(\mathbf{r}_{i_1} \mathbf{r}_{i_1}^{\top}-N^{-1} \boldsymbol{\Phi}\right) \mathbb{E}_j\left(\mathbf{Q}_{i_1j}^{-1}\left(z_1\right)\check{\mathbf{Q}}_{i_1j}^{-1}\left(z_2\right)\right)\\
&\qquad\qquad\qquad\quad\times 
\tr\mathbf{H}_n\left(z_1\right)\left(\mathbf{r}_{i_2} \mathbf{r}_{i_2}^{\top}-N^{-1} \boldsymbol{\Phi}\right) \mathbb{E}_j\left(\mathbf{Q}_{i_2j}^{-1}\left(z_1\right)\check{\mathbf{Q}}_{i_2 j}^{-1}\left(z_2\right)\right),
\end{align*}
where $\check{\bQ}_{i_2j}$ is defined similarly as $\bQ_{i_2j}$ by $(\br_1,\ldots,\br_{j-1},\check{\br}_{j+1},\ldots,\check{\br}_{p})$
and where $\check{\br}_{j+1},\ldots,\check{\br}_{p}$ are i.i.d. copies of $\br_{j+1},\ldots,\br_{p}$.

When $i_1=i_2$, with Lemma 5 in \cite{gao2017high}, the term in the above expression is bounded by 
$$
\left|b_1\left(z_1\right)\right|^2 \sum_{i_1<j} \mathbb{E}
\left|\tr\mathbf{H}_n\left(z_1\right)\left(\mathbf{r}_{i_1} \mathbf{r}_{i_1}^{\top}-N^{-1} \boldsymbol{\Phi}\right) \mathbb{E}_j\left(\mathbf{Q}_{i_1 j}^{-1}\left(z_1\right)\check{\mathbf{Q}}_{i_1 j}^{-1}\left(z_2\right)\right)\right|^2\leq Kpc_n^{-1}n^{-1}=O(1).
$$
For $i_1\neq i_2<j$, define
$$
\beta_{i_1i_2j}(z)=\frac{1}{\sqrt{c_N}+\mathbf{r}_{i_2}^{\top} \mathbf{Q}_{i_1i_2 j}^{-1}(z) \mathbf{r}_{i_2}},~ 
\mathbf{Q}_{i_1i_2 j}(z)=\bQ(z)-\sqrt{\frac{N}{p}}(\mathbf{r}_{i_1} \mathbf{r}_{i_1}^{\top}+\mathbf{r}_{i_2} \mathbf{r}_{i_2}^{\top}+\mathbf{r}_j \mathbf{r}_j^{\top}),
$$
and similarly define $\check{\beta}_{i_1i_2 j}$ and $\check{\bQ}_{i_1i_2 j}(z)$. Then we have
\begin{align*}
\left|b_1(z_1)\right|^2\sum_{i_1\neq i_2<j}&\E\tr\mathbf{H}_n\left(z_1\right)\left(\mathbf{r}_{i_1} \mathbf{r}_{i_1}^{\top}-N^{-1} \boldsymbol{\Phi}\right) \mathbb{E}_j\left(\mathbf{Q}_{i_1j}^{-1}\left(z_1\right)\check{\mathbf{Q}}_{i_1j}^{-1}\left(z_2\right)\right)\\
&\times 
\tr\mathbf{H}_n\left(z_1\right)\left(\mathbf{r}_{i_2} \mathbf{r}_{i_2}^{\top}-N^{-1} \boldsymbol{\Phi}\right) \mathbb{E}_j\left(\mathbf{Q}_{i_2j}^{-1}\left(z_1\right)\check{\mathbf{Q}}_{i_2 j}^{-1}\left(z_2\right)\right)=S_1+S_2+S_3,
\end{align*}
where \begin{align*}
S_1&=-\sum_{i_1\neq i_2<j}\E\tr\mathbf{H}_n\left(z_1\right)\left(\mathbf{r}_{i_1} \mathbf{r}_{i_1}^{\top}-N^{-1} \boldsymbol{\Phi}\right) \mathbb{E}_j\left(
\beta_{i_1i_2j}(z_1)\bQ_{i_1i_2j}^{-1}(z_1)\br_{i_2}\br_{i_2}^{\top}\bQ_{i_1i_2j}^{-1}(z_1)
\check{\mathbf{Q}}_{i_1j}^{-1}\left(z_2\right)\right)\\
&\qquad\qquad\qquad\times 
\left|b_1(z_1)\right|^2
\tr\mathbf{H}_n\left(z_1\right)\left(\mathbf{r}_{i_2} \mathbf{r}_{i_2}^{\top}-N^{-1} \boldsymbol{\Phi}\right) \mathbb{E}_j\left(\mathbf{Q}_{i_2j}^{-1}\left(z_1\right)\check{\mathbf{Q}}_{i_2 j}^{-1}\left(z_2\right)\right),\\
S_2&=-\sum_{i_1\neq i_2<j}\E\tr\mathbf{H}_n\left(z_1\right)\left(\mathbf{r}_{i_1} \mathbf{r}_{i_1}^{\top}-N^{-1} \boldsymbol{\Phi}\right) \mathbb{E}_j\left(
\bQ_{i_1i_2j}^{-1}(z_1)
\check{\beta}_{i_1i_2j}(z_2)\check{\mathbf{Q}}_{i_1i_2j}^{-1}\left(z_2\right)\br_{i_2}\br_{i_2}^{\top}\check{\mathbf{Q}}_{i_1i_2j}^{-1}\left(z_2\right)
\right)\\
&\qquad\qquad\qquad\times \left|b_1(z_1)\right|^2
\tr\mathbf{H}_n\left(z_1\right)\left(\mathbf{r}_{i_2} \mathbf{r}_{i_2}^{\top}-N^{-1} \boldsymbol{\Phi}\right) \mathbb{E}_j\left(\mathbf{Q}_{i_2j}^{-1}\left(z_1\right)\check{\mathbf{Q}}_{i_2 j}^{-1}\left(z_2\right)\right),\\
S_3&=-\left|b_1(z_1)\right|^2\sum_{i_1\neq i_2<j}\E\tr\mathbf{H}_n\left(z_1\right)\left(\mathbf{r}_{i_1} \mathbf{r}_{i_1}^{\top}-N^{-1} \boldsymbol{\Phi}\right) \mathbb{E}_j\left(
\bQ_{i_1i_2j}^{-1}(z_1)
\check{\mathbf{Q}}_{i_1i_2j}^{-1}(z_2)
\right)\\
&\qquad\qquad\qquad\times 
\tr\mathbf{H}_n\left(z_1\right)\left(\mathbf{r}_{i_2} \mathbf{r}_{i_2}^{\top}-N^{-1} \boldsymbol{\Phi}\right) \mathbb{E}_j\left(\beta_{i_2 i_1j}\left(z_1\right) \mathbf{Q}_{i_1 i_2 j}^{-1}\left(z_1\right) \mathbf{r}_{i_1} \mathbf{r}_{i_1}^{\top}\mathbf{Q}_{i_1 i_2j}^{-1}\left(z_1\right)\check{\mathbf{Q}}_{i_2 j}^{-1}\left(z_2\right)\right).
\end{align*}
Spilt $$
S_1=S_1^{(1)}+S_1^{(2)},~S_1^{(2)}=S_1^{(21)}+S_1^{(22)},~S_1^{(22)}=S_1^{(221)}+S_1^{(222)},
$$
where
\begin{align*}
    S_1^{(1)}&=\sum_{i_1\neq i_2<j}\E\tr\mathbf{H}_n\left(z_1\right)\left(\mathbf{r}_{i_1} \mathbf{r}_{i_1}^{\top}-N^{-1} \boldsymbol{\Phi}\right)\\
    &\qquad\times 
    \mathbb{E}_j\left(
\beta_{i_1i_2j}(z_1)\bQ_{i_1i_2j}^{-1}(z_1)\br_{i_2}\br_{i_2}^{\top}\bQ_{i_1i_2j}^{-1}(z_1)
\check{\beta}_{i_1 i_2 j}\left(z_2\right) \check{\mathbf{Q}}_{i_1 i_2 j}^{-1}\left(z_2\right) \mathbf{r}_{i_2} \mathbf{r}_{i_2}^{\top} \check{\mathbf{Q}}_{i_1 i_2 j}^{-1}\left(z_2\right)\right)\\
&\qquad\times 
\left|b_1(z_1)\right|^2
\tr\mathbf{H}_n\left(z_1\right)\left(\mathbf{r}_{i_2} \mathbf{r}_{i_2}^{\top}-N^{-1} \boldsymbol{\Phi}\right) \mathbb{E}_j\left(\mathbf{Q}_{i_2j}^{-1}\left(z_1\right)\check{\mathbf{Q}}_{i_2 j}^{-1}\left(z_2\right)\right),\\
S_1^{(2)}&=-\sum_{i_1\neq i_2<j}\E\tr\mathbf{H}_n\left(z_1\right)\left(\mathbf{r}_{i_1} \mathbf{r}_{i_1}^{\top}-N^{-1} \boldsymbol{\Phi}\right) \mathbb{E}_j\left(
\beta_{i_1i_2j}(z_1)\bQ_{i_1i_2j}^{-1}(z_1)\br_{i_2}\br_{i_2}^{\top}\bQ_{i_1i_2j}^{-1}(z_1)
\check{\mathbf{Q}}_{i_1i_2j}^{-1}\left(z_2\right)\right)\\
&\qquad\qquad\qquad\times 
\left|b_1(z_1)\right|^2
\tr\mathbf{H}_n\left(z_1\right)\left(\mathbf{r}_{i_2} \mathbf{r}_{i_2}^{\top}-N^{-1} \boldsymbol{\Phi}\right) \mathbb{E}_j\left(\mathbf{Q}_{i_2j}^{-1}\left(z_1\right)\check{\mathbf{Q}}_{i_2 j}^{-1}\left(z_2\right)\right),\\
S_1^{(21)}&=\sum_{i_1\neq i_2<j}\E\tr\mathbf{H}_n\left(z_1\right)\left(\mathbf{r}_{i_1} \mathbf{r}_{i_1}^{\top}-N^{-1} \boldsymbol{\Phi}\right) \mathbb{E}_j\left(
\beta_{i_1i_2j}(z_1)\bQ_{i_1i_2j}^{-1}(z_1)\br_{i_2}\br_{i_2}^{\top}\bQ_{i_1i_2j}^{-1}(z_1)
\check{\mathbf{Q}}_{i_1i_2j}^{-1}\left(z_2\right)\right)\\
&\quad\times 
\left|b_1(z_1)\right|^2
\tr\mathbf{H}_n\left(z_1\right)\left(\mathbf{r}_{i_2} \mathbf{r}_{i_2}^{\top}-N^{-1} \boldsymbol{\Phi}\right) \mathbb{E}_j\left(\beta_{i_2i_1j}\left(z_1\right) \mathbf{Q}_{i_2 i_1 j}^{-1}\left(z_1\right) \mathbf{r}_{i_1} \mathbf{r}_{i_1}^{\top}\mathbf{Q}_{i_2 i_1 j}^{-1}\left(z_1\right)\check{\mathbf{Q}}_{i_2 j}^{-1}\left(z_2\right)\right),\\
S_1^{(22)}&=-\sum_{i_1\neq i_2<j}\E\tr\mathbf{H}_n\left(z_1\right)\left(\mathbf{r}_{i_1} \mathbf{r}_{i_1}^{\top}-N^{-1} \boldsymbol{\Phi}\right) \mathbb{E}_j\left(
\beta_{i_1i_2j}(z_1)\bQ_{i_1i_2j}^{-1}(z_1)\br_{i_2}\br_{i_2}^{\top}\bQ_{i_1i_2j}^{-1}(z_1)
\check{\mathbf{Q}}_{i_1i_2j}^{-1}\left(z_2\right)\right)\\
&\qquad\qquad\qquad\times 
\left|b_1(z_1)\right|^2
\tr\mathbf{H}_n\left(z_1\right)\left(\mathbf{r}_{i_2} \mathbf{r}_{i_2}^{\top}-N^{-1} \boldsymbol{\Phi}\right) \mathbb{E}_j\left(\mathbf{Q}_{i_2i_1j}^{-1}\left(z_1\right)\check{\mathbf{Q}}_{i_2j}^{-1}\left(z_2\right)\right),\\
S_1^{(221)}&=\sum_{i_1\neq i_2<j}\E\tr\mathbf{H}_n\left(z_1\right)\left(\mathbf{r}_{i_1} \mathbf{r}_{i_1}^{\top}-N^{-1} \boldsymbol{\Phi}\right) \mathbb{E}_j\left(
\beta_{i_1i_2j}(z_1)\bQ_{i_1i_2j}^{-1}(z_1)\br_{i_2}\br_{i_2}^{\top}\bQ_{i_1i_2j}^{-1}(z_1)
\check{\mathbf{Q}}_{i_1i_2j}^{-1}\left(z_2\right)\right)\\
&\times 
\left|b_1(z_1)\right|^2
\tr\mathbf{H}_n\left(z_1\right)\left(\mathbf{r}_{i_2} \mathbf{r}_{i_2}^{\top}-N^{-1} \boldsymbol{\Phi}\right) \mathbb{E}_j\left(\mathbf{Q}_{i_2i_1j}^{-1}\left(z_1\right)
\check{\beta}_{i_2 i_1 j}\left(z_1\right) \check{\mathbf{Q}}_{i_2 i_1 j}^{-1}\left(z_1\right) \mathbf{r}_{i_1} \mathbf{r}_{i_1}^{\top} \check{\mathbf{Q}}_{i_2 i_1 j}^{-1}\left(z_1\right)
\right),\\
S_1^{(222)}&=-\sum_{i_1\neq i_2<j}\E\tr\mathbf{H}_n\left(z_1\right)\left(\mathbf{r}_{i_1} \mathbf{r}_{i_1}^{\top}-N^{-1} \boldsymbol{\Phi}\right) \mathbb{E}_j\left(
\beta_{i_1i_2j}(z_1)\bQ_{i_1i_2j}^{-1}(z_1)\br_{i_2}\br_{i_2}^{\top}\bQ_{i_1i_2j}^{-1}(z_1)
\check{\mathbf{Q}}_{i_1i_2j}^{-1}\left(z_2\right)\right)\\
&\qquad\qquad\qquad\times 
\left|b_1(z_1)\right|^2
\tr\mathbf{H}_n\left(z_1\right)\left(\mathbf{r}_{i_2} \mathbf{r}_{i_2}^{\top}-N^{-1} \boldsymbol{\Phi}\right) \mathbb{E}_j\left(\mathbf{Q}_{i_2i_1j}^{-1}\left(z_1\right)\check{\mathbf{Q}}_{i_2i_1j}^{-1}\left(z_2\right)\right).
\end{align*}
With \eqref{estimate}, we have 
$$
S_1^{(1)}=O(n),~S_1^{(21)}=O(n),~S_1^{(221)}=O(n),~S_1^{(222)}=0,$$
which gives us $S_1=O(n)$. Similarly, we can show $S_2=O(n),~S_3=O(n)$. Hence, we obtain \eqref{A3-o1}.

For $A_1(z_1,z_2)$, we have $$
\mathbb{E}\left|A_1\left(z_1, z_2\right)+\frac{j-1}{n^2} b_1\left(z_2\right) \operatorname{tr}\left(\mathbb{E}_j\left(\mathbf{Q}_j^{-1}\left(z_1\right)\right) \mathbf{Q}_j^{-1}\left(z_2\right)\right) \operatorname{tr}\left(\mathbf{Q}_j^{-1}\left(z_2\right) \mathbf{H}_n\left(z_1\right)\right)\right| \leq K p^{1 / 2},
$$
from which we obtain 
\begin{align*}
&\operatorname{tr}\left(\mathbb{E}_j\left(\mathbf{Q}_j^{-1}\left(z_1\right)\right) \mathbf{Q}_j^{-1}\left(z_2\right)\right)
=-\operatorname{tr}\left(\mathbf{H}_n\left(z_1\right) \mathbf{Q}_j^{-1}\left(z_2\right)\right)\nonumber\\
&-\frac{j-1}{n^2}b_1\left(z_1\right) b_1\left(z_2\right) \operatorname{tr}\left(\mathbb{E}_j\left(\mathbf{Q}_j^{-1}\left(z_1\right)\right) \mathbf{Q}_j^{-1}\left(z_2\right)\right)\operatorname{tr}\left(\mathbf{Q}_j^{-1}\left(z_2\right) \mathbf{H}_n\left(z_1\right)\right)
+A_4(z_1,z_2),
\end{align*}
where $\E|A_4(z_1,z_2)|\leq Kn^{1/2}$. By the similar strategy as the proof of \eqref{A3-o1}, we have $$
\E|\E_j\tr b_1(z_2)\bA(z_2)\bH_n(z_1)|\leq \sqrt{\E|\E_j\tr b_1(z_2)\bA(z_2)\bH_n(z_1)|^2}\leq Kn^{-1/2},
$$
from which we obtain
\begin{align*}
&\operatorname{tr}\left(\mathbb{E}_j\left(\mathbf{Q}_j^{-1}\left(z_1\right)\right) \mathbf{Q}_j^{-1}\left(z_2\right)\right)
=\operatorname{tr}\left(\mathbf{H}_n\left(z_1\right) \mathbf{H}_n\left(z_2\right)\right)\\&+\frac{j-1}{n^2}b_1\left(z_1\right) b_1\left(z_2\right) \operatorname{tr}\left(\mathbb{E}_j\left(\mathbf{Q}_j^{-1}\left(z_1\right)\right) \mathbf{Q}_j^{-1}\left(z_2\right)\right)
\operatorname{tr}\left(\mathbf{H}_n\left(z_2\right) \mathbf{H}_n\left(z_1\right)\right)+A_5(z_1,z_2),
\end{align*}
where $\E|A_5(z_1,z_2)|\leq Kn^{1/2}+Knp^{-1/2}$. Define $g_n(z)=\sqrt{c_N}\E\beta_1(z)$. 
Similar as (B.40)-(B.41) in \cite{gao2017high}, we then have
$$
g_n(z)=\sqrt{c_N}\E\beta_1(z)=-\left(c_N+\sqrt{c_N} z\right) \frac{1}{c_n} \E\left[\frac{1}{\sqrt{c_N}} s_n(z)+\frac{1-c_n}{c_N+\sqrt{c_N} z}\right]
$$
And with \eqref{estimate},we have
\begin{align}\label{gn}
|\sqrt{c_N}b_n(z)-g_n(z)|=\sqrt{c_N}|b_n(z)-\E\beta_1(z)|\leq K p^{-1/2},~
|b_n(z)-b_1(z)|\leq Kn^{1/2}p^{-3/2},
\end{align}
which implies
\begin{align*}
&\operatorname{tr}\left(\mathbb{E}_j\left(\mathbf{Q}_j^{-1}\left(z_1\right)\right) \mathbf{Q}_j^{-1}\left(z_2\right)\right)
=\operatorname{tr}\left(\mathbf{H}_n\left(z_1\right) \mathbf{H}_n\left(z_2\right)\right)\\&+\frac{j-1}{n^2}
\frac{g_n\left(z_1\right) g_n\left(z_2\right)}{c_N}\operatorname{tr}\left(\mathbb{E}_j\left(\mathbf{Q}_j^{-1}\left(z_1\right)\right) \mathbf{Q}_j^{-1}\left(z_2\right)\right)
\operatorname{tr}\left(\mathbf{H}_n\left(z_2\right) \mathbf{H}_n\left(z_1\right)\right)+A_6(z_1,z_2),
\end{align*}
where $\mathbb{E}\left|A_6\left(z_1, z_2\right)\right| \leq K n^{1 / 2}+K n p^{-1 / 2}$. Recall $\mathbf{H}_n\left(z\right)=\left(\sqrt{c_N}+z-\frac{p-1}{N} b_1\left(z\right)\right)^{-1} \mathbf{I}_n$ and let $$d_n(z_1,z_2)=\frac{1}{n}\tr \bH_n(z_1)\bH_n(z_2),~
a_n(z_1,z_2)=g_n(z_1)g_n(z_2)d_n(z_1,z_2).$$ 
Since $c_nb_n\left(z_1\right) b_n\left(z_2\right) d_n\left(z_1, z_2\right)/a_n\left(z_1, z_2\right)\to 1$, $\mathcal{J}$ can be written as 
$$
\mathcal{J}=\frac{1}{p}a_n(z_1,z_2)\sum_{j=1}^{p}\frac{1}{1-\frac{j-1}{p}a_n(z_1,z_2)}+A_7(z_1,z_2),
$$
where $\mathbb{E}\left|A_7\left(z_1, z_2\right)\right| \leq K n^{-1 / 2}$. Note that the limit of $a_n\left(z_1, z_2\right)$ is $a\left(z_1, z_2\right)=\frac{1}{\left(s(z_1)+z_1\right)\left(s(z_2)+z_2\right)}$. Thus the limit of $\frac{\partial^2}{\partial z_2 \partial z_1} \mathcal{J}$ in probability is
\begin{align*}
\frac{\partial^2}{\partial z_2 \partial z_1} \int_0^{a\left(z_1, z_2\right)} \frac{1}{1-z} d z&=\frac{\partial}{\partial z_2}\left(\frac{\partial a\left(z_1, z_2\right) / \partial z_1}{1-a\left(z_1, z_2\right)}\right)=\frac{s'(z_1)s'(z_2)}{\left[s(z_1)-s(z_2)\right]^2}-\frac{1}{(z_1-z_2)^2}\\
&=\frac{s^2\left(z_1\right) s^2\left(z_2\right)}{\left[s^2\left(z_1\right)-1\right]\left[s^2\left(z_2\right)-1\right]\left[s\left(z_1\right) s\left(z_2\right)-1\right]^2}.
\end{align*}
%Similar as the argument for $\mathcal{J}$ and Lemma 9 in %\cite{gao2017high}, we have 
%$$
%\frac{\partial^2}{\partial z_2 \partial z_1} %J_4\xrightarrow{\text { i.p. }}\frac{|\E %x_{11}^2|^2 s^2\left(z_1\right) %s^2\left(z_2\right)}%%{\left[s^2\left(z_1\right)-1\right]\left[s^2\left(z%_2\right)-1\right]\left[|\E x_{11}^2|^2 %s\left(z_1\right) s\left(z_2\right)-1\right]^2}.
%$$
Thus we complete the proof of Lemma \ref{le-J}.

%under the real case. %Similarly, for the complex case we obtain %$M_n^{(2)}(z)=\eqref{complex-%theta}+\frac{\sqrt{c_N}}%{z\left(\sqrt{c_N}+z\right)}$.

	\begin{supplement}
		\stitle{Supplementary Material of ``On  eigenvalues of a renormalized sample correlation matrix"}
		\sdescription{This supplementary document contains the proofs of 
		Lemmas \ref{le-max-eig}, \ref{step2},\ref{step3},\ref{le-supp}.
}\end{supplement}

%%  MR numbers will be added by VTeX.                      %%
%%                                                         %%
%%  Use \cite{...} to cite references in text.             %%
%%                                                         %%
%%%%%%%%%%%%%%%%%%%%%%%%%%%%%%%%%%%%%%%%%%%%%%%%%%%%%%%%%%%%%

%% if your bibliography is in bibtex format, uncomment commands:
\bibliographystyle{imsart-nameyear} % Style BST file (imsart-number.bst or imsart-nameyear.bst)
\bibliography{reference}       % Bibliography file (usually '*.bib')

%% or include bibliography directly:
% \begin{thebibliography}{}
% \bibitem[\protect\citeauthoryear{???}{???}]{b1}
% \end{thebibliography}

\end{document}